\documentclass[final,onefignum,onetabnum]{siamart220329}

\usepackage{braket,amsfonts}

\usepackage{array}
\usepackage{bm}
\usepackage{amsmath,amssymb}
\usepackage{mathtools}
\usepackage{stmaryrd}
\usepackage[caption=false]{subfig}
\usepackage{pgfplots}

\newsiamthm{claim}{Claim}
\newsiamremark{remark}{Remark}
\newsiamremark{hypothesis}{Hypothesis}
\crefname{hypothesis}{Hypothesis}{Hypotheses}

\usepackage[noend]{algpseudocode}

\usepackage{footnote}
\usepackage{tikz}
\usepackage{mathdots}
\usepackage{yhmath}
\usepackage{pifont}
\usetikzlibrary{fadings}
\usetikzlibrary{patterns}
\usetikzlibrary{shadows.blur}
\usetikzlibrary{shapes}
\usetikzlibrary{matrix}
\usetikzlibrary{shapes.geometric, arrows}
\usetikzlibrary{matrix,positioning,decorations.pathreplacing}
\tikzstyle{io} = [rectangle, rounded corners, 
minimum width=1.2cm, 
minimum height=1.3cm,
text centered, 
text width = 3.5cm,
draw=black, 
fill=blue!5]
\tikzstyle{ior} = [rectangle, rounded corners, 
minimum width=1.2cm, 
minimum height=1.3cm,
text centered, 
text width = 3cm,
draw=black, 
fill= red!5]
\tikzstyle{iorr} = [rectangle, rounded corners, 
minimum width=1.2cm, 
minimum height=1.3cm,
text centered, 
text width = 3cm,
draw=black, 
fill= red!20]
\tikzstyle{check} = [rectangle, rounded corners, 
minimum width=0.7cm, 
minimum height=0.7cm,
text centered, 
text width = 0.7cm,
draw=black, 
fill=green!10]
\tikzstyle{rnk} = [rectangle, rounded corners, 
minimum width=1cm, 
minimum height=1cm,
text centered, 
text width = 3.5cm,
draw=black, 
fill=red!5]
\tikzstyle{rnkch} = [rectangle, rounded corners, 
minimum width=1cm, 
minimum height=1cm,
text centered, 
text width = 3cm,
draw=black, 
fill=blue!5]
\tikzstyle{arrow} = [thick,->,>=stealth]

\usepackage{graphicx,epstopdf}

\Crefname{ALC@unique}{Line}{Lines}
\allowdisplaybreaks

\usepackage{amsopn}

\DeclareMathOperator{\rk}{rank}

\newcommand{\R}{\mathbb{R}}

\newcommand{\bigO}{\mathcal{O}}

\newcommand{\lowrank}[2]{\llbracket {#1} \rrbracket_{#2}}
\newcommand{\norm}[1]{\left\lVert#1\right\rVert}

\newcommand{\acur}{AdaCUR}
\newcommand{\qcur}{FastAdaCUR}

\usepackage[normalem]{ulem}

\newcommand{\ignore}[1]{}

\usepackage{xspace}
\usepackage{bold-extra}
\usepackage[most]{tcolorbox}

\colorlet{texcscolor}{blue!50!black}
\colorlet{texemcolor}{red!70!black}
\colorlet{texpreamble}{red!70!black}
\colorlet{codebackground}{black!25!white!25}

\lstdefinestyle{siamlatex}{%
  style=tcblatex,
  texcsstyle=*\color{texcscolor},
  texcsstyle=[2]\color{texemcolor},
  keywordstyle=[2]\color{texemcolor},
  moretexcs={cref,Cref,maketitle,mathcal,text,headers,email,url},
}

\tcbset{%
  colframe=black!75!white!75,
  coltitle=white,
  colback=codebackground, 
  colbacklower=white, 
  fonttitle=\bfseries,
  arc=0pt,outer arc=0pt,
  top=1pt,bottom=1pt,left=1mm,right=1mm,middle=1mm,boxsep=1mm,
  leftrule=0.3mm,rightrule=0.3mm,toprule=0.3mm,bottomrule=0.3mm,
  listing options={style=siamlatex}
}

\newtcblisting[use counter=example]{example}[2][]{%
  title={Example~\thetcbcounter: #2},#1}

\newtcbinputlisting[use counter=example]{\examplefile}[3][]{%
  title={Example~\thetcbcounter: #2},listing file={#3},#1}

\DeclareTotalTCBox{\code}{ v O{} }
{ 
  fontupper=\ttfamily\color{black},
  nobeforeafter,
  tcbox raise base,
  colback=codebackground,colframe=white,
  top=0pt,bottom=0pt,left=0mm,right=0mm,
  leftrule=0pt,rightrule=0pt,toprule=0mm,bottomrule=0mm,
  boxsep=0.5mm,
  #2}{#1}

\patchcmd\newpage{\vfil}{}{}{}
\begin{tcbverbatimwrite}{tmp_\jobname_header.tex}
\title{Low-rank approximation of parameter-dependent matrices via CUR decomposition
\thanks{Date: \today 
\funding{TP is supported by the Heilbronn Institute for Mathematical Research. YN is supported by EPSRC grants EP/Y010086/1 and EP/Y030990/1.}
}
}

\author{Taejun Park\thanks{Mathematical Institute, University of Oxford, Oxford, OX2 6GG, UK, (\email{park@maths.ox.ac.uk}, \email{nakatsukasa@maths.ox.ac.uk}).}
\and Yuji Nakatsukasa\footnotemark[2] }

\headers{Low-rank approx of parameter-dependent matrices via CUR}{Taejun Park and Yuji Nakatsukasa}
\end{tcbverbatimwrite}
\title{Low-rank approximation of parameter-dependent matrices via CUR decomposition
\thanks{Date: \today
\funding{TP is supported by the Heilbronn Institute for Mathematical Research. YN is supported by EPSRC grants EP/Y010086/1 and EP/Y030990/1.}
}
}

\author{Taejun Park\thanks{Mathematical Institute, University of Oxford, Oxford, OX2 6GG, UK, (\email{park@maths.ox.ac.uk}, \email{nakatsukasa@maths.ox.ac.uk}).}
\and Yuji Nakatsukasa\footnotemark[2] }

\headers{Low-rank approx of parameter-dependent matrices via CUR}{Taejun Park and Yuji Nakatsukasa}

\ifpdf
\hypersetup{ pdftitle={Low-rank approximation of parameter-dependent matrices via CUR} }
\fi

\begin{document}
\maketitle

\begin{tcbverbatimwrite}{tmp_\jobname_abstract.tex}
\begin{abstract}
A low-rank approximation of a parameter-dependent matrix $A(t)$ is an important task in the computational sciences appearing for example in dynamical systems and compression of a series of images. In this work, we introduce \acur{}, an efficient algorithm for computing a low-rank approximation of parameter-dependent matrices via CUR decompositions. The key idea for this algorithm is that for nearby parameter values, the column and row indices for the CUR decomposition can often be reused. \acur{} is rank-adaptive, provides error control, and has complexity that compares favorably against existing methods. A faster algorithm which we call \qcur{} that prioritizes speed over accuracy is also given, which is rank-adaptive and has complexity which is at most linear in the number of rows or columns, but without error control.
\end{abstract}


\end{tcbverbatimwrite}
\begin{abstract}
A low-rank approximation of a parameter-dependent matrix $A(t)$ is an important task in the computational sciences appearing for example in dynamical systems and compression of a series of images. In this work, we introduce \acur{}, an efficient algorithm for computing a low-rank approximation of parameter-dependent matrices via CUR decompositions. The key idea for this algorithm is that for nearby parameter values, the column and row indices for the CUR decomposition can often be reused. \acur{} is rank-adaptive, provides error control, and has complexity that compares favorably against existing methods. A faster algorithm which we call \qcur{} that prioritizes speed over accuracy is also given, which is rank-adaptive and has complexity which is at most linear in the number of rows or columns, but without error control.
\end{abstract}




\section{Introduction}
The task of finding a low-rank approximation to a matrix is ubiquitous in the computational sciences \cite{lowrank_udelltownsend}. In large scale problems, low-rank approximation provides an efficient way to store and process the matrix. In this work, we study the low-rank approximation of a parameter-dependent matrix $A(t) \in \R^{m\times n}$ with $m\geq n$ for a finite number of parameter values $t\in \R$ in some compact domain $D \subset \R$. This task has appeared in several applications, for example, in compression of a series of images \cite{NonnenmacherLubich2008}, dynamical systems \cite{KochLubich2007} and Gaussian random fields \cite{KressnerLatzMasseiUllmann2020}.

When $A(t) \equiv A$ is a constant, the truncated singular value decomposition (TSVD) attains the best low-rank approximation in any unitarily invariant norm for a fixed target rank $r$ \cite[\S~7.4.9]{hornjohn}. However, the cost of computing the truncated SVD becomes infeasible for large matrices. The situation is exacerbated when $A(t)$ varies with parameter $t$, requiring the evaluation of the truncated SVD at every parameter value $t$ of interest. Consequently, numerous alternative approaches have been explored in recent years such as dynamical low-rank approximation \cite{KochLubich2007,NonnenmacherLubich2008}, randomized SVD and generalized Nystr\"om method \cite{KressnerLam2023} and the CUR decomposition \cite{DonelloPalkarNaderiDelReyFernandezBabaee2023}. 

In this work, we use the CUR decomposition \cite{GoreinovTyrtyshinikovZamarashkin1997,MahoneyDrineas2009} to find a low-rank approximation to a parameter-dependent matrix $A(t)$. The CUR decomposition of a matrix $A$ is a low-rank approximation that takes the form
\begin{equation}
  A\approx A(:,J) A(I,J)^\dagger A(I,:) = CU^\dagger R =: A_{IJ}
\end{equation} in MATLAB notation. 
Here, $C = A(:,J)$ is a subset of the columns of $A$ with $J$ being the column indices, $R = A(I,:)$ is a subset of the rows of $A$ with $I$ being the row indices and $U = A(I,J)$ is the intersection of $C$ and $R$.\footnote{There are other choices for the CUR decomposition. In particular, the choice $A\approx CC^\dagger A R^\dagger R$ minimizes the Frobenius norm error given $C$ and $R$. We choose $A = CU^\dagger R$ in this paper as it is faster to compute.} While the CUR decomposition may be considered suboptimal in comparison to the truncated SVD, it possesses many favorable properties. Firstly, given the row and column indices, they are extremely efficient to compute and store as they do not require reading the entire matrix $A$. Additionally, as the low-rank factors are the subsets of the original matrix $A$, they inherit certain properties of $A$ such as sparsity and non-negativity. Moreover, these factors assist with data interpretation by revealing the important rows and columns. For further insight and theory on the CUR decomposition, see, for example, \cite{BoutsidisWoodruff2017,HammHuang2020,MahoneyDrineas2009}.

For any set of row indices $I$ with $|I| = r$ and column indices $J$ with $|J| = r$ where $r$ is the target rank, the CUR decomposition satisfies the following error bound \cite{pn24cur}
\begin{equation} \label{eq:CAbound}
    \norm{A-CU^\dagger R}_F \leq \norm{Q_C(I,:)^\dagger}_2 \norm{Q_X(J,:)^{-1}}_2 \norm{A-AX^\dagger X}_F
\end{equation} where $X\in \R^{r\times n}$ is any row space approximator of $A$ and $Q_X \in \R^{n\times r}$ and $Q_C \in \R^{m\times r}$ are orthonormal matrices spanning the columns of $X^T$ and $C = A(:,J)$ respectively.\footnote{This bound can also be stated in terms of a column approximator and a subset of the rows of $A$.} The first two terms on the right-hand side of the inequality $\norm{Q_C(I,:)^\dagger}_2$ and $\norm{Q_X(J,:)^{-1}}_2$ mostly govern the accuracy of the CUR decomposition, because the row space approximator $X$ is arbitrary and can be chosen such that $\norm{A-AX^\dagger X}_F$ is quasi-optimal. Therefore, it is important to get a good set of row and column indices that control the first two terms. There are many existing algorithms that lead to a good set of indices such as leverage score sampling \cite{DrineasMahoneyMuthukrishnan2008,MahoneyDrineas2009}, DPP sampling \cite{DerezinskiMahoney2021}, volume sampling \cite{CortinovisKressner2020,DeshpandeRademacherVempalaWang2006} and pivoting strategies \cite{DongMartinsson2023,DrmacGugercin2016,matcomp,GuEisenstat1996,SorensenEmbree2016,TrefethenBau,VoroninMartinsson2017}, which all attempt to minimize the first two terms. Notably, there exists a set of row and column indices for which the CUR decomposition satisfies the following near-optimal guarantee \cite{ZamarashkinOsinsky2018}
\begin{equation*}
    \norm{A - CU^\dagger R}_F \leq (1+r) \norm{A-\lowrank{A}{r}}_F
\end{equation*} where $\lowrank{A}{r}$ is the best rank-$r$ approximation to $A$. 
This means that little is lost by requiring our low-rank approximation to be in CUR form, as long as the indices are chosen appropriately.
Polynomial-time algorithms for constructing such indices can be found in \cite{CortinovisKressner2020,Osinsky2023}. The bound \eqref{eq:CAbound} can be improved by oversampling the rows, that is, by obtaining more row indices such that $I$ has more indices than the target rank. The topic of oversampling has been explored in \cite{AndersonDuMahoneyMelgaardWuGu2015,GidisuHochstenbach2022,pn24cur,PeherstorferDrmacGugercin2020,ZimmermannWillcox2016}. Notably, \cite{pn24cur} describes how to compute the CUR decomposition in a numerically stable way and shows that oversampling can help improve the stability (in addition to accuracy) of the CUR decomposition.

Now, having defined what the CUR decomposition is, the objective of this work is the following:
\begin{objective}
    Let $A(t)$ be a parameter-dependent matrix. Then given a set of parameter values $t_1,t_2,...,t_q$, and a tolerance $\epsilon$, devise an algorithm that approximately achieves
    \begin{equation}
    \norm{A(t_i) - \big(A(t_i)(:,J_i)\big)\left(A(t_i)(I_i,J_i)^\dagger \right)\left(A(t_i)(I_i,:)\right)}_F \leq \epsilon \norm{A(t_i)}_F
\end{equation} for each $i$, where $I_i$ and $J_i$ are the row and column indices corresponding to $t_i$.
\end{objective}

For parameter-dependent matrices $A(t)$, we can naively recompute the row and the column indices for each parameter value of interest. However, this approach can be inefficient and wasteful, as the computed indices for one parameter value are likely to contain information about nearby parameter values. The reason is because a matrix may only undergo slight deviations when the parameter value changes. More concretely, let $A(t)$ be a parameter-dependent matrix and $I$ and $J$ be a set of row and column indices respectively. Then by setting $X = V_r(t)^T$ where $V_r(t)$ contains the $r$-dominant right singular vectors of $A(t)$, \eqref{eq:CAbound} changes to
\begin{equation} \label{eq:paramCAbound}
    \norm{A(t)-C(t) U(t)^\dagger R(t)}_F \leq \norm{Q_{C(t)}(I,:)^\dagger}_2 \norm{Q_{V_r(t)}(J,:)^{-1}}_2 \norm{A(t)-\lowrank{A(t)}{r}}_F
\end{equation} where $\lowrank{A(t)}{r}$ is the best rank-$r$ approximation to $A(t)$. Now the rightmost term in \eqref{eq:paramCAbound} is optimal for any parameter value $t$, so we compare the first two terms. Suppose we have two different parameter values, say $t_1$ and $t_2$. Then, if both $\norm{Q_{C(t_1)}(I,:)^\dagger}_2$ $\norm{Q_{V_r(t_1)}(J,:)^{-1}}_2$ and $\norm{Q_{C(t_2)}(I,:)^\dagger}_2 \norm{Q_{V_r(t_2)}(J,:)^{-1}}_2$ remain small with the same set of indices $I$ and $J$, we can use $I$ and $J$ to approximate both $A(t_1)$ and $A(t_2)$.\footnote{In our work, we allow $\norm{A(t_1)-A(t_2)}_2$ to be large, say $\bigO(1)$, and therefore the perturbation bounds for the CUR decomposition in \cite{HammHuang2021} is not enough to explain why the same indices can be reused for nearby parameter values.}

More specifically, define the following set of parameter-dependent matrices:
\begin{equation}
    \mathcal{M}_{IJ}^{(r)}(\delta):= \left\{A(t) \in \R^{m\times n}: \norm{Q_{A(t)(:,J)}(I,:)^\dagger}_2 \norm{Q_{V(t)}(J,:)^{-1}}_2 \leq \delta \text{ for all }t \in D \right\}
\end{equation} for some domain $D \subset \R$ and sets of indices $I$ and $J$ with $|I| = |J| = r$. Then, for all $A(t) \in \mathcal{M}_{IJ}^{(r)}(\delta)$, 
\begin{equation} 
    \norm{A(t)-C(t) U(t)^\dagger R(t)}_F \leq \delta \norm{A(t)-\lowrank{A(t)}{r}}_F
\end{equation} holds for all $t\in D$. Therefore if $A(t) \in \mathcal{M}_{IJ}^{(r)}(\delta)$ for sufficiently small $\delta \geq 1$ then the index sets $I$ and $J$ provide a good rank-$r$ CUR approximation to $A(t)$ for all $t\in D$. An example of a class of parameter-dependent matrices belonging to $\mathcal{M}_{IJ}^{(r)}(\delta)$ are incoherent rank-$r$ parameter-dependent matrices. In general, it is difficult to precisely determine whether a parameter-dependent matrix belongs to $\mathcal{M}_{IJ}^{(r)}(\delta)$ for some (practical) starting indices $I$ and $J$ or not, but intuitively, any parameter-dependent matrix that is sufficiently low-rank and has singular vectors that change gradually should belong to $\mathcal{M}_{IJ}^{(r)}(\delta)$ for some $I,J,\delta$ and $r$. We do not explore the technicalities of this further in this work.

Nevertheless, in practice, for nearby parameter values, the identical set of indices $I$ and $J$ often capture the important rows and columns of both $A(t_1)$ and $A(t_2)$; see Sections \ref{sec:method} and \ref{sec:numill}. This motivates us to reuse the indices for various parameter values, potentially resulting in significant savings, often up to a factor $r$ in complexity where $r$ is the target rank. With this key idea in mind, the main goal of this paper is to devise a rank-adaptive certified\footnote{We use the term \emph{certified} to denote approximation methods that are designed to control the approximation error.} algorithm that reuses the indices as much as possible until the error becomes too large, at which point we recompute the indices.

To achieve this goal efficiently and reliably, we rely on a variety of efficient tools, many of which are borrowed from the randomized linear algebra literature such as randomized rank estimation \cite{AndoniNguyen2013,MeierNakatsukasa2024} and pivoting on a random sketch \cite{DongMartinsson2023,Duersch_Gu_SIREV,VoroninMartinsson2017}; see Section \ref{sec:prelim}. We assume throughout that for each parameter $t$, $A(t)$ has decaying singular values so that a low-rank approximation is possible. This is clearly a required assumption for low-rank approximation techniques to be effective.

In the next section, we provide brief outlines of our algorithms, \acur{} and \qcur{}, review existing methods for approximating parameter-dependent matrices, and summarize our contributions. This is followed by an overview of the various techniques from (randomized) numerical linear algebra that we use in our algorithms. The next two sections cover \acur{} and \qcur{} in detail, as well as numerical experiments highlighting their performance. Finally, we conclude and discuss topics for further research.

\subsection{Brief outline of \acur{} and \qcur{}}
Before introducing the two algorithms, \acur{} and \qcur{}, in detail in Section \ref{sec:method}, we first provide a brief outline of their core ideas.

\acur{} is aimed at efficiently computing an accurate low-rank approximation of parameter-dependent matrices. A high-level overview of \acur{} is given in Figure~\ref{adacurOverview}.

\begin{figure} \label{adacurOverview}
\footnotesize
\begin{tcolorbox}[colback=white,colframe=black,boxsep=2pt]{\centering \underline{\textbf{\acur{}}} 
\smallskip
 \begin{enumerate}
        \item[\textbf{Input}]: Parameter-dependent matrix $A(t)\in \R^{m\times n}$, parameters $t_1,...,t_q$, error tol. $\epsilon$.
        \smallskip
        \item[\textbf{Output}]: CUR factors $C(t_j), U(t_j), R(t_j)$ such that $A(t)\approx C(t_j)U(t_j)^\dagger R(t_j)$.
        \medskip
        \item Compute initial set of indices for $A(t_1)$.
        \smallskip
        \item[\ding{98}] For $j = 2,3,...,q$,\vspace{0.3cm}  \\ 
        \resizebox{4.42in}{!}{
        \begin{tikzpicture}[node distance=1.5cm,scale = 0.5] 
\hspace{-1cm}
\node (relerr) [io] {Compute rel. error for $A(t_j)$ using prev. indices};
\node (fin1) [check, below of=relerr,yshift = -0.2cm] {\ding{52}};
\node (minmod) [ior, right of=relerr,xshift = 4.2cm] {Make minor modifications to the indices};
\node (recomp) [iorr, right of=minmod, xshift = 3.8cm] {Recompute indices from scratch};
\node (fin2) [check, below of = minmod, yshift = -0.2cm] {\ding{52}};

\draw [arrow] (relerr) -- node[anchor=west] {Rel. Err. $\leq \epsilon$}(fin1);
\draw [arrow] (relerr) -- node[anchor=south] {Rel. Err. $> \epsilon$}(minmod);
\draw [arrow] (minmod) -- node[anchor=west] {Rel. Err. $\leq \epsilon$}(fin2);
\draw [arrow] (minmod) -- node[anchor=south] {Rel. Err. $> \epsilon$}(recomp);
\end{tikzpicture}
}
\end{enumerate}}
\end{tcolorbox}
\caption{An overview of \acur{}}
\end{figure}

\acur{} begins by computing an initial set of indices for $A(t_1)$. For subsequent parameter values $t_j$, \acur{} first computes the relative error for $A(t_j)$ using the previous sets of indices. If the error is smaller than the tolerance $\epsilon$, the algorithm retains the previous sets of indices and proceeds to the next parameter value. If the error exceeds the tolerance, minor modifications are made to the indices, and the algorithm checks whether the tolerance is satisfied. If the modified indices meet the tolerance, they are used; otherwise, the indices are recomputed from scratch. The precise details of \acur{} are outlined in Section \ref{subsec:mainalg}.

\qcur{} is aimed at achieving even greater speed at the cost of possible reduction in accuracy. Instead of computing the error which can be expensive, \qcur{} keeps a small set of extra indices that helps with rank adaptivity. A high-level overview of \qcur{} is given in Figure~\ref{fastadacurOverview}.

\begin{figure} \label{fastadacurOverview}
\footnotesize
\begin{tcolorbox}[colback=white,colframe=black,boxsep=2pt]
{\centering \underline{\textbf{\qcur{}}}
\smallskip
 \begin{enumerate}
        \item[\textbf{Input}]: Parameter-dependent matrix $A(t)\in \R^{m\times n}$, parameters $t_1,...,t_q$, rank tol. $\epsilon$
        \smallskip
        \item[\textbf{Output}]: CUR factors $C(t_j), U(t_j), R(t_j)$ such that $A(t)\approx C(t_j)U(t_j)^\dagger R(t_j)$.
        \medskip
        \item Compute initial set of indices $I,J$ and a small set of extra indices $I_s, J_s$ for $A(t_1)$.
        \item[\ding{98}] For $j = 2,3,...,q$,\vspace{0.3cm}  \\ 
        \resizebox{3.5in}{!}{
        \hspace{1cm}
        \begin{tikzpicture}[node distance=1.5cm,scale = 0.5]
\hspace{-1cm}
\node (coreMat) [rnkch] {Compute $\epsilon$-rank of $A(I\cup I_s,J\cup J_s)$};
\node (rankincr) [rnk, right of=coreMat,xshift = 3.8cm] {Add indices to $I,J$ from $I_s,J_s$ and add more indices to $I_s,J_s$};
\node (rankdec) [rnk, below of=rankincr,yshift = 0.3cm] {Remove indices from $I,J$};

\draw [arrow] (coreMat) -- node[anchor=north] {$\epsilon$-rank decr.}(rankdec);
\draw [arrow] (coreMat) -- node[anchor=south] {$\epsilon$-rank incr.}(rankincr);
\end{tikzpicture}
}
\end{enumerate}}
\end{tcolorbox}
\caption{An overview of \qcur{}}
\end{figure}

\qcur{} begins by computing an initial sets of indices $I,J$ and a small set of additional indices $I_s,J_s$ for $A(t_1)$. For subsequent parameter values $t_j$, \qcur{} calculates the $\epsilon$-rank of the core matrix $A(I\cup I_s,J\cup J_s)$ using the previous sets of indices to determine whether the $\epsilon$-rank exceeds the size of the index sets $I$ and $J$. Importantly, for \qcur{}, the entire matrix is not accessed beyond the first parameter value, which makes the algorithm efficient. If the $\epsilon$-rank has increased, indices are added to $I$ and $J$ from $I_s$ and $J_s$. Conversely, if the $\epsilon$-rank has decreased, indices are removed from $I$ and $J$. After this adjustment, the algorithm proceeds to the next parameter value. The details of \qcur{} are provided in Section \ref{subsec:cheapalg}. Unlike \acur{}, \qcur{} does not have error control mechanism, making it vulnerable to adversarial examples; see Section \ref{subsec:adversarial}.

\paragraph{Existing methods} \label{para:existing}
There have been several different approaches for finding a low-rank approximation to parameter-dependent matrices $A(t)$. We describe four different classes of methods. First, dynamical low-rank approximation is a differential-equation-based approach where the low-rank factors are obtained from projecting the time derivative of $A(t)$, $\dot{A}(t)$, onto the tangent space of the smooth manifold consisting of fixed-rank matrices; see \cite{Lubich2014}. A number of variations from the original dynamical low-rank approximation \cite{KochLubich2007} have been proposed to deal with various issues such as the stiffness introduced by the curvature of the manifold \cite{CerutiLubich2022,KieriLubichWalach2016,LubichOseledets2014} and rank-adaptivity \cite{CerutiKuschLubich2022,HauckSchnake2023,HesthavenPagliantiniRipamonti2022,HochbruckNeherSchrammer2023}. The complexity of this approach is typically $\bigO(r_i T_{A(t_i)})$ for the parameter $t_i$, where $r_i$ is the target rank for $A(t_i)$ and $T_{A(t_i)}$ is the cost of matrix-vector product with $A(t_i)$ or $A(t_i)^T$. Secondly, Kressner and Lam \cite{KressnerLam2023} use randomized SVD \cite{HalkoMartinssonTropp2011} and generalized Nystr\"om \cite{Nakatsukasa2020,TroppYurtseverUdellCevher2017}, both based on random projections, to find low-rank approximations of parameter-dependent matrices $A(t)$. They focus on matrices that admit an affine linear decomposition with respect to $t$. Notably, they use the same random embedding for all parameter values rather than generating a new random embedding for each parameter value. The complexity is also $\bigO(r_i T_{A(t_i)})$ for the parameter value $t_i$. Next, Donello et al. \cite{DonelloPalkarNaderiDelReyFernandezBabaee2023} use the CUR decomposition to find a low-rank approximation to parameter-dependent matrices. At each iteration, a new set of column and row indices for the CUR decomposition are computed by applying the discrete empirical interpolation method (DEIM) \cite{ChaturantabutSorenson2010,DrmacGugercin2016,SorensenEmbree2016} to the singular vectors of the CUR decomposition from the previous iteration. While this approach benefits from not viewing the entire matrix, it may suffer from the same adversarial example as the fast algorithm, \qcur{} discussed in this paper; see Algorithm \ref{alg:cheap} and Section \ref{subsec:adversarial}. They allow oversampling of column and/or row indices to improve the accuracy of the CUR decomposition and propose rank-adaptivity where the rank is either increased or decreased by $1$ if the smallest singular value of the current iterate is larger or smaller, respectively, than some threshold. Although their algorithm is based on the CUR decomposition, the low-rank factors in their approach are represented in the form of an SVD. This SVD is computed from the CUR decomposition at each iteration, given that orthonormal factors are necessary in DEIM. The complexity for this algorithm is $\bigO((m+n)r_i^2)$ for the parameter value $t_i$. Lastly, in the special case when $A(t)$ is symmetric positive definite, a parametric variant of adaptive cross approximation has been used in \cite{KressnerLatzMasseiUllmann2020}.

\paragraph{Contributions}
The contribution of this paper lies in the design of efficient rank-adaptive algorithms for low-rank approximations of parameter-dependent matrices. The algorithms are based on the CUR decomposition of a matrix. We first demonstrate that the same set of row and column indices, or making slight modifications to the indices, can yield effective low-rank approximations for nearby parameter values. This observation forms the basis of our efficient algorithms, which we describe in Section \ref{sec:method}. We present two distinct algorithms: a rank-adaptive algorithm with error control, which we call \acur{} (Adaptive CUR), and a faster rank-adaptive algorithm that sacrifices error control for speed, which we name \qcur{} (Fast Adaptive CUR); see Algorithms \ref{alg:main} and \ref{alg:cheap} respectively. 

\acur{} has the worst-case time complexity of $\bigO\left(r_iT_{A(t_i)} + (m+n)r_i^2\right)$ for the parameter value $t_i$ where $r_i$ is the target rank for $A(t_i)$ and $T_{A(t_i)}$ denotes the cost of matrix-vector product with $A(t_i)$ or $A(t_i)^T$. However, in practice, the algorithm frequently runs with the best-case time complexity of $\bigO(T_{A(t_i)}+(m+n)r_i)$ for the parameter value $t_i$; see Sections \ref{subsec:mainalg} and \ref{sec:numill}. This is competitive with existing methods such as \cite{CerutiLubich2022,KressnerLam2023,LubichOseledets2014}, which run with complexity $\bigO(r_iT_{A(t_i)})$. Notably, our algorithm's rank-adaptive nature and built-in error control offer distinct advantages over many existing algorithms, which often lack one or both of these features.

\qcur{}, aside from the initial phase of computing the initial indices, runs linearly in $m$ and $n$ in the worse-case. The best-case complexity is $\bigO(r_i^3)$ for the parameter value $t_i$, which makes the algorithm remarkably fast. While this algorithm is also rank-adaptive, it lacks rigorous error control since its priority is efficiency over accuracy. However, the algorithm has the advantage of needing only $\mathcal{O}(r)$ rows and columns for each parameter value; i.e., there is no need to view the full matrix. This feature is particularly attractive when the entries are expensive to compute. In experiments, we notice that the error may grow as we iterate through the parameter for difficult problems. Here, difficult problems refer to problems where $A(t)$ undergoes frequent large rank changes or one that changes rapidly such that the previous indices carry little to no information for the next. Nevertheless, we frequently observe that the algorithm performs well on easier problems; see Section \ref{sec:numill}.

\paragraph{Notation}
Throughout, we use $\norm{\cdot}_2$ for the spectral norm or the vector-$\ell_2$ norm and $\norm{\cdot}_F$ for the Frobenius norm. We use dagger $^\dagger$ to denote the pseudoinverse of a matrix and $\lowrank{A}{r}$ to denote the best rank-$r$ approximation to $A$ in any unitarily invariant norm, i.e., the approximation derived from truncated SVD~\cite[\S~7.4.9]{hornjohn}. We use $a_1 \lesssim a_2$ to denote that $a_1 \leq C a_2$ for some constant $C>0$. Unless specified otherwise, $\sigma_i(A)$ denotes the $i$th largest singular value of the matrix $A$. We use MATLAB style notation for matrices and vectors. For example, for the $k$th to $(k+j)$th columns of a matrix $A$ we write $A(:,k:k+j)$. Lastly, we use $|I|$ to denote the length of the vector or the cardinality of the set $I$, $I_1-I_2$ to be the set difference between $I_1$ and $I_2$, and define $[n] :=\{1,2,...,n\}$. 

\section{Preliminaries} \label{sec:prelim}
In this section, we discuss the tools needed for our proposed methods introduced in Section \ref{sec:method}. Many of the methods described here are in the randomized numerical linear algebra literature. For an in-depth review, we refer to \cite{HalkoMartinssonTropp2011,MartinssonTropp2020,Woodruff2014}. First, in Section \ref{subsec:randemb}, we discuss the core ingredient in many randomized numerical linear algebra algorithms, which are random embeddings. Following this, in Section \ref{subsec:pivot}, we review pivoting on a random sketch, which efficiently computes the indices for the CUR decomposition, and in Section \ref{subsec:rankest}, we review a fast randomized algorithm for computing the numerical rank of a matrix. Lastly, in Section \ref{subsec:normest}, we discuss an efficient method for computing the Frobenius norm of a matrix, which is used for error estimation for \acur{} in Section \ref{sec:method}.

\subsection{Random embeddings} \label{subsec:randemb}
Let $A\in \R^{m\times n}$ be a matrix of rank $r \leq \min\{m,n\}$ (typically $r \ll \min\{m,n\}$). Then $\Gamma \in \R^{s\times m}$ $(r\leq s \leq \min\{m,n\})$ is a subspace embedding for the span of $A$ with distortion $\epsilon \in (0,1)$ if
\begin{equation} \label{eq:subemb}
    (1-\epsilon) \norm{Ax}_2 \leq \norm{\Gamma Ax}_2 \leq (1+\epsilon) \norm{Ax}_2
\end{equation} for every $x\in \R^n$. Therefore, $\Gamma$ is a linear map which preserves the $2$-norm of every vector in a given subspace. A random embedding is a subspace embedding drawn at random that satisfies \eqref{eq:subemb} for all $A$ with high probability. A typical application of random embeddings is in \emph{matrix sketching}, a technique for dimensionality reduction. Matrix sketching compresses the original matrix into a smaller-sized matrix while retaining much of its original information. For example, for a matrix $A\in \R^{m\times n}$ of rank $r \ll \min\{m,n\}$, $\Gamma A \in \R^{s\times n}$ where $r\leq s\ll \min\{m,n\}$ is a (row) sketch of the matrix $A$. Here, the integer $s$ is called the sketch size.

There are a few important examples of random embeddings such as Gaussian embeddings \cite{HalkoMartinssonTropp2011,MartinssonTropp2020}, subsampled randomized trigonometric transforms (SRTTs) \cite{BoutsidisGittens2013,Tropp2011} and sparse embeddings \cite{ClarksonWoodruff2017,cohen16,NelsonNguyen2013}. We focus on Gaussian embeddings in this work. A Gaussian embedding $\Gamma \in \R^{s\times m}$ has i.i.d. entries $\Gamma_{ij} \sim \mathcal{N}(0,1/s)$. Here, $s = \Omega(r/\epsilon^2)$ for theoretical guarantees, but $s = r+\Omega(1)$ works well in practice. The cost of applying a Gaussian embedding to a matrix $A\in \R^{m\times n}$ is $\bigO(s T_A)$ where $T_A$ is the cost of matrix-vector multiply with $A$ or $A^T$. Other random embeddings provide a similar theoretical guarantee but generally require a larger sketch size, say $s = \Omega(r\log r/\epsilon^2)$. However, they are usually cheaper to apply; for example, SRTTs cost $\bigO(mn\log s)$ \cite{WoolfeLibertyRokhlinTygert2008}.

\subsection{Pivoting on a random sketch} \label{subsec:pivot}
Pivoting on a random sketch \cite{DongMartinsson2023,Duersch_Gu_SIREV,VoroninMartinsson2017} is an attractive and efficient method for finding a spanning set of columns and/or rows of a matrix. This will be the core part for computing the CUR decomposition. This method has two main steps: sketch and pivot. Let $A\in \R^{m\times n}$ be a matrix and $r$ be the target rank. Then the basic idea works as follows for column indices:
\begin{enumerate}
    \item \textit{Sketch}: Draw a random embedding $\Gamma\in \R^{r\times m}$ and form $X = \Gamma A \in \R^{r\times n}$.\footnote{For robustness, oversampling is recommended, that is, we draw a random embedding with the sketch size larger than $r$, say $2r$, in step $1$ and obtain $r$ pivots in step $2$. See \cite{DongChenMartinssonPearce2024,DongMartinsson2023} for a discussion.}
    \item \textit{Pivot}: Perform a pivoting scheme, e.g., column pivoted QR (CPQR) or partially pivoted LU (LUPP) on $X$. Collect the chosen pivot indices.
\end{enumerate}
The sketching step entails compressing the original matrix down to a smaller-sized matrix using a random embedding. The sketch $X = \Gamma A$ forms a good approximation to the row space of $A$ \cite{HalkoMartinssonTropp2011}. The pivoting step involves applying a pivoting scheme on the smaller-sized matrix $X$, which reduces the cost compared to applying the pivoting scheme directly on $A$. There are many pivoting schemes such as column pivoted QR (CPQR), partially pivoted LU (LUPP) or those with strong theoretical guarantees such as Gu-Eisenstat's strong rank-revealing QR (sRRQR) \cite{GuEisenstat1996} or BSS sampling \cite{BatsonSpielmanSrivastava2012,BoutsidisDrineasMagdon-Ismail2014}. See \cite{DongMartinsson2023} for a comparison. In this work, we use a version of Algorithm $1$ from \cite{DongMartinsson2023}, which we state below in Algorithm \ref{alg:randpivot}.

\begin{algorithm}[!ht]\footnotesize
  \caption{Pivoting on a random sketch}
  \label{alg:randpivot}
  \begin{algorithmic}[1]
\Require{$A \in \R^{m\times n}$, target rank $r$ (typically $r \ll \min\{m,n\}$)}
\Ensure{Column indices $J$ and row indices $I$ with $|I| = |J| = r$}
\vspace{0.5pc}
\Statex \texttt{function} $[I,J] = \mathtt{Rand\_ Pivot}(A,r)$
\State Draw a Gaussian embedding $\Gamma\in \R^{r\times m}$.
\State Set $X = \Gamma A \in \R^{r\times n}$, a row sketch of $A$
\State Apply CPQR on $X$. Let $J$ be the $r$ column pivots.
\State Apply CPQR on $A(:,J)^T$. Let $I$ be the $r$ row pivots.
\end{algorithmic}
\end{algorithm}
Algorithm \ref{alg:randpivot} selects the column indices first by applying CPQR on the sketch $X = \Gamma A$, and then selects the row indices by applying CPQR on the chosen columns $A(:,J)^T$.\footnote{The fourth step in Algorithm \ref{alg:randpivot}, which involves applying CPQR on the selected columns instead of, for example, a column sketch of $A$, can be important for the stability and accuracy of the CUR decomposition \cite{pn24cur}.} The complexity of Algorithm \ref{alg:randpivot} is $\bigO(rT_A + (m+n)r^2)$, where $rT_A$ is the cost of forming the Gaussian row sketch $X$ in line $2$ (see Section \ref{subsec:randemb}) and $\bigO(nr^2)$ and $\bigO(mr^2)$ are the cost of applying CPQR on $X$ (line $3$) and $A(:,J)^T$ (line $4$) respectively. 

\subsection{Randomized rank estimation} \label{subsec:rankest}
Randomized rank estimation \cite{AndoniNguyen2013,MeierNakatsukasa2024} is an efficient tool for finding the $\epsilon$-rank of a matrix $A\in \R^{m\times n}$. As we often require the target rank as input, for example in Algorithm \ref{alg:randpivot}, we would like an efficient method for computing the $\epsilon$-rank of a matrix. This method also relies on sketching but, unlike pivoting on a random sketch, it involves sketching $\Gamma_1,\Gamma_2$ from both sides of the matrix. The idea is to approximate the $\epsilon$-rank of $A$ using
\begin{equation} \label{eq:rankest}
    \hat{r}_\epsilon = \rk_{\epsilon}(\Gamma_1 A \Gamma_2)
\end{equation} where $\rk_{\epsilon}(B)$ is the $\epsilon$-rank of the matrix $B$ (i.e., the number of singular values larger than $\epsilon$) and $\Gamma_1 A \Gamma_2 \in \R^{s\times 2s}$ with $s \gtrsim \rk_{\epsilon}(A)$. Here, $s \gtrsim \rk_{\epsilon}(A)$ is achieved by gradually increasing $s$ by appending to the sketches until $\sigma_{\min}(\Gamma_1 A \Gamma_2) <\epsilon$; see \cite{MeierNakatsukasa2024} for details. The overall cost of the algorithm is $\bigO(s T_A+ns\log s + s^3)$ where $\bigO(s T_A)$ is the cost of forming the Gaussian sketch, $\bigO(ns \log s )$ is the cost of forming the SRTT sketch and $\bigO(s^3)$ is the cost of computing the singular values of $\Gamma_1 A \Gamma_2$.

\paragraph{Combining pivoting on a random sketch and randomized rank estimation}
Since we require the target rank as input for pivoting on a random sketch (Algorithm \ref{alg:randpivot}), randomized rank estimation and pivoting on a random sketch can be used together. Both algorithms need to form the row sketch $X = \Gamma A$ as part of their algorithm. Therefore, we can input the row sketch $X$ obtained from randomized rank estimation into pivoting on a random sketch. This avoids the need to reform the row sketch making randomized rank estimation almost free when used with pivoting on a random sketch. The overall algorithm is outlined in Algorithm \ref{alg:pivotrankest}, giving the overall complexity of $\bigO(\hat{r}_\epsilon T_A + (m+n)\hat{r}_\epsilon^2)$. By combining the two algorithms into one, the number of matrix-vector products with $A$ is halved, which is often the dominant cost.

\begin{algorithm}[!ht]\footnotesize
  \caption{Pivoting on a random sketch using randomized rank estimation}
  \label{alg:pivotrankest}
  \begin{algorithmic}[1]
\Require{$A \in \R^{m\times n}$, tolerance $\epsilon$}
\Ensure{Row indices $I$ and Column indices $J$}
\vspace{0.5pc}
\Statex \texttt{function} $[I,J] = \mathtt{Rand\_Pivot\_RankEst}(A,\epsilon)$
\State $[\hat{r}_\epsilon,X] = $ Estimate the $\epsilon$-rank rank using \eqref{eq:rankest} \Comment{$X = \Gamma_1 A$ is the row sketch from \eqref{eq:rankest}}
\State $[I,J] = \mathtt{Rand\_Pivot}(A,X,\hat{r}_\epsilon )$ \Comment{Take the row sketch $X$ as input in Algorithm \ref{alg:randpivot}}
\end{algorithmic}
\end{algorithm}

\subsection{Randomized norm estimation} \label{subsec:normest}
In order to monitor how well our algorithm performs, we need an efficient way to estimate the norm of the error. Randomized norm estimation \cite{GrattonTitley-Peloquin2018,MartinssonTropp2020} is an efficient method for finding an approximation to a norm of a matrix. In this work, we focus on the Frobenius norm as it is a natural choice for low-rank approximation and easier to estimate. Let $A\in \R^{m\times n}$ be a matrix and $\hat{A}_r \in \R^{m\times n}$ be a rank-$r$ approximation to $A$. Then we would like to estimate $\norm{A-\hat{A}_r}_F$. As before, we can sketch to get an approximation to $A-\hat{A}_r$ and then compute the Frobenius norm. It turns out that the sketch size of $O(1)$, say $5$, suffices. More specifically, we have the following theorem from \cite{GrattonTitley-Peloquin2018}.

\begin{theorem}{\cite[Theorem 3.1]{GrattonTitley-Peloquin2018}}
    Let $A\in \R^{m\times n}$ be a matrix, $\Gamma\in \R^{s\times m}$ be a Gaussian matrix with i.i.d. entries $\Gamma_{ij} \sim \mathcal{N}(0,1)$ and set $\rho:= \norm{A}_F^2/\norm{A}_2^2$. For any $\tau > 1$ and $s\leq m$,
    \begin{equation}
        \Pr\left(\frac{\norm{A}_F}{\tau}<\frac{\norm{\Gamma A}_F}{\sqrt{s}} \leq \tau \norm{A}_F\right) \geq 1-\exp\left(-\frac{s\rho}{2}(\tau-1)^2\right) - \exp\left(-\frac{s\rho}{4}\frac{(\tau^2-1)^2}{\tau^4}\right)
    \end{equation}
\end{theorem}

Setting $\norm{A}_F = 5\norm{A}_2$, $\tau = 2$ and $s = 5$, the above theorem tells us
\begin{equation}
 \Pr\left(\frac{\norm{A}_F}{2}<\frac{\norm{\Gamma A}_F}{\sqrt{5}} \leq 2 \norm{A}_F\right) \geq 1-2.33\cdot 10^{-8},   
\end{equation} i.e. with a small number $s$ of Gaussian vectors, we can well-approximate the Frobenius norm of a matrix. The algorithm specialized to CUR decomposition is presented in Algorithm \ref{alg:normest}.

\begin{algorithm}[!ht]\footnotesize
  \caption{Randomized norm estimation for CUR decomposition}
  \label{alg:normest}
  \begin{algorithmic}[1]
\Require{$A \in \R^{m\times n}$, $I$ row indices, $J$ column indices, sample size $s$ (rec. $s = 5$) }
\Ensure{$E$, an estimate for $\norm{A-A(:,J)A(I,J)^\dagger A(I,:)}_F$}
\vspace{0.5pc}
\Statex \texttt{function} $E = \mathtt{Norm\_ Est}(A,I,J,s)$
\State Draw a Gaussian embedding $\Gamma \in \R^{s \times m}$
\State Set $X = \Gamma A \in \R^{s\times n}$
\State Set $\hat{X} = \Gamma  A(:,J)A(I,J)^\dagger A(I,:)\in \R^{s\times n}$
\State $E = \norm{X-\hat{X}}_F$
\end{algorithmic}
\end{algorithm}
The complexity of Algorithm \ref{alg:normest} is $\bigO(sT_A + (m+n)rs) = \bigO(T_A+(m+n)r)$ where $r = \max\{|I|,|J|\}$ since $s = \bigO(1)$. This method is particularly useful when $A$ can only be accessed through matrix-vector multiply.

In addition to approximating the error using randomized norm estimation, we obtain, as a by-product, a small sketch of the low-rank residual matrix since
\begin{equation}
    X - \hat{X} = \Gamma(A - A(:,J)A(I,J)^\dagger A(I,:)).
\end{equation} This small sketch can be used to obtain an additional set of row and column indices by using pivoting on it, similarly to Algorithm \ref{alg:randpivot}. Specifically, we apply pivoting on the sketched residual $X-\hat{X}$ to get an extra set of $s$ column indices, denoted $J_s$. Next, we apply pivoting on the chosen columns of the residual matrix,
\begin{equation*}
    A(:,J_s) - A(:,J)A(I,J)^\dagger A(I,J_s)
\end{equation*} to extract an extra set of $s$ row indices $I_s$. It is important to note that $I_s$ and $J_s$ will not include indices already chosen in $I$ and $J$ because the residual matrix $A-A(:,J)A(I,J)^\dagger A(I,:)$ is zero in the rows and columns corresponding to $I$ and $J$. The procedure for obtaining $s$ additional indices from the sketched residual is outlined in Algorithm \ref{alg:pivotRes}.

\begin{algorithm}[!ht]\footnotesize
  \caption{Pivoting on the residual matrix}
  \label{alg:pivotRes}
  \begin{algorithmic}[1]
\Require{$X_R \in \R^{s\times n}$ sketch of the residual matrix, $I$ row indices, $J$ column indices}
\Ensure{$I_s$ and $J_s$, extra sets of indices}
\vspace{0.5pc}
\Statex \texttt{function} $[I_s,J_s] = \mathtt{Pivot\_Residual}(X_R,I,J)$
\State Apply CPQR on $X_R$. Let $J_s$ be the $s$ column pivots,
\State Set $C_R = A(:,J_s) - A(:,J)A(I,J)^\dagger A(I,J_s) \in \R^{m\times s}$
\State Apply CPQR on $C_R^T$. Let $I_s$ be the $s$ row pivots.
\end{algorithmic}
\end{algorithm}

The complexity of Algorithm \ref{alg:pivotRes} is $\bigO((m+n)s^2+mrs+r^3) = \bigO(n+mr+r^3)$, where $r = \max\{|I|,|J|\}$ since $s = \bigO(1)$. The dominant cost lies in the computation of $A(:,J)A(I,J)^\dagger A(I,J_s)$, which requires $\bigO(mr+r^3)$ operations. By adding these additional indices, we can improve the accuracy of the CUR decomposition. Algorithm \ref{alg:pivotRes} will be applied later in \acur{} to refine the approximation.

\paragraph{Multiple error estimation}
Algorithm \ref{alg:normest} can easily be used to estimate multiple low-rank approximations of $A$. Suppose we are given two sets of indices $I_1$ and $J_1$, and $I_2$ and $J_2$. Then once the row sketch of $A$, $X = \Gamma A$ has been computed, we only need to compute the row sketches of the low-rank approximations $\hat{X}_1 = \Gamma A_{I_1 J_1}$ and $\hat{X}_2 = \Gamma A_{I_2 J_2}$ to approximate the error. Here, $A_{IJ} = A(:,J)A(I,J)^\dagger A(I,:)$. Since forming $X$ is typically the dominant cost in Algorithm \ref{alg:normest}, this allows us to efficiently test multiple low-rank approximations.

\section{Proposed method} \label{sec:method}
Let us first restate the problem. Let $A(t) \in \R^{m\times n}$ be a parameter-dependent matrix with $t\in D$ where $D \subseteq \R$ is a compact domain and $t_1,t_2,...,t_q \in D$ be a finite number of distinct sample points ordered in some way, e.g., $t_1 < t_2 < \cdots <t_q$ for $t\in \R$. Then find a low-rank approximation of $A(t_1),A(t_2),...,A(t_q)$.\footnote{Although our algorithms are presented for a finite number of sample points, they are likely to work on a continuum, e.g. by choosing the indices for the closest point $t_i$ for each $t\in[t_1,t_q]$.} 

The goal of this work is to devise an efficient algorithm that reuses the indices to their fullest extent in the CUR decomposition as we iterate through the parameter values. To see how we can use the same indices for different parameter values, recall \eqref{eq:CAbound} with oversampling from \cite{pn24cur}:
\begin{equation} \label{eq:CAboundOS}
    \norm{A-A_{I\cup I_0,J}}_F \leq \norm{Q_C(I \cup I_0,:)^\dagger}_2 \norm{V_r(J,:)^{-1}}_2 \norm{A-\lowrank{A}{r}}_F.
\end{equation}
Here, $I_0$ is an extra set of row indices resulting from oversampling that are distinct from $I$.

The bound \eqref{eq:CAboundOS} leads us to two observations.
First, as described earlier in the introduction, if $A(t) \in \mathcal{M}_{IJ}^{(r)}(\delta)$, then 
\begin{equation*}
    \norm{A-A_{I\cup I_0,J}}_F \leq \delta \norm{A- \lowrank{A}{r}}_F.
\end{equation*} Therefore, if $\delta$ is sufficiently small, say $\bigO(10)$, the sets of indices $I\cup I_0$ and $J$ yield a good CUR approximation. However, this does not imply that a high value of $\delta$ necessarily results in a poor CUR approximation. Let us illustrate this with a numerical example. We use the synthetic example from \cite{CerutiLubich2022} given by 
    \begin{equation} \label{eq:synprob}
    A(t) = e^{tW_1}e^t D e^{tW_2}, t\in [0,1]
\end{equation} where $D\in \R^{n\times n}$ is diagonal with entries $d_{jj} = 2^{-j}$ and $W_1,W_2\in \R^{n\times n}$ are randomly generated skew-symmetric matrices. The singular values of $A(t)$ are $e^t 2^{-j}$ for $j= 1,2,...,n$. We take $n = 500$ and take $101$ equispaced points in the interval $[0,1]$ for $t$. The matrix exponentials were computed using the $\mathtt{expm}$ command in MATLAB.

\begin{figure}[!ht]
\vspace{-0.5cm}
\subfloat[\centering Rel. error plot: CUR and trunc. SVD]{\label{subfig:SameP}\includegraphics[scale = 0.33]{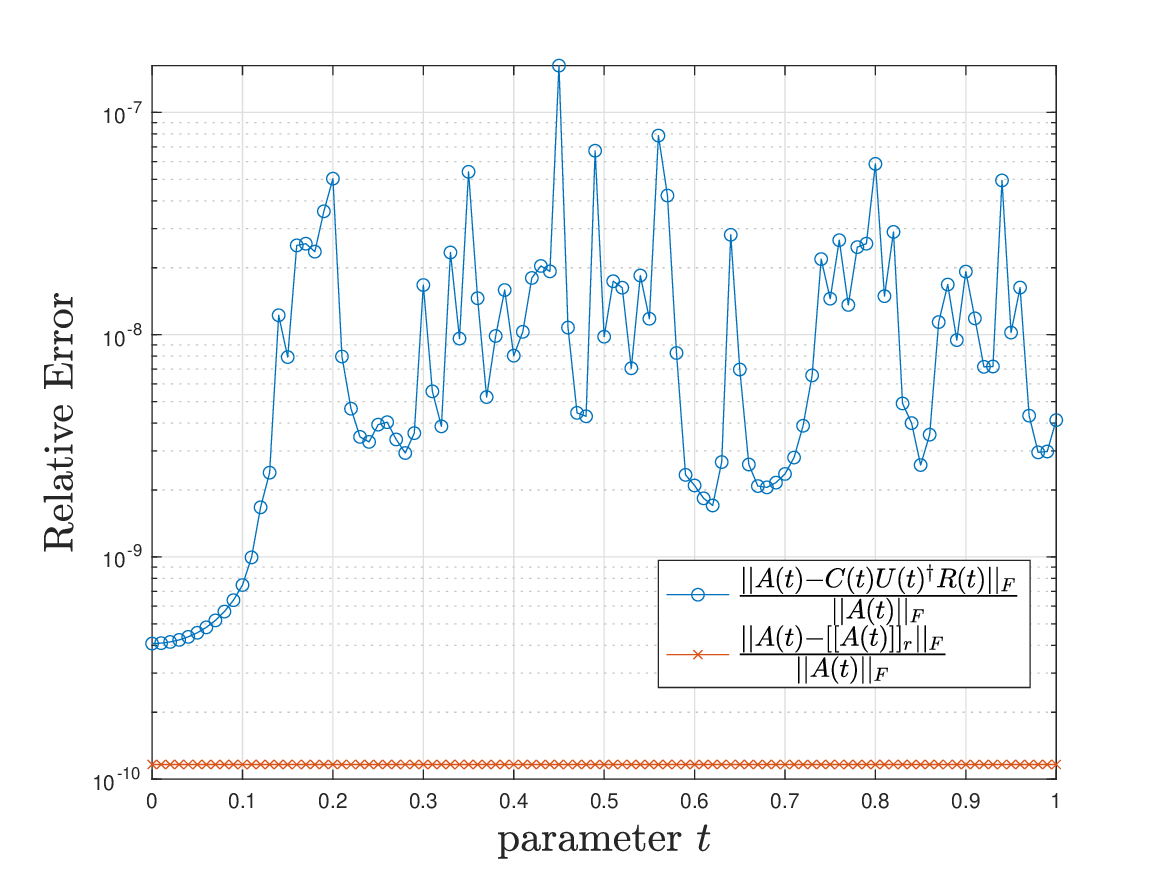}} \hspace{-0.5cm}
\subfloat[\centering Gap between theory and practice of the CUR bound in \eqref{eq:CAboundOS}.]{\label{subfig:SameC}\includegraphics[scale = 0.33]{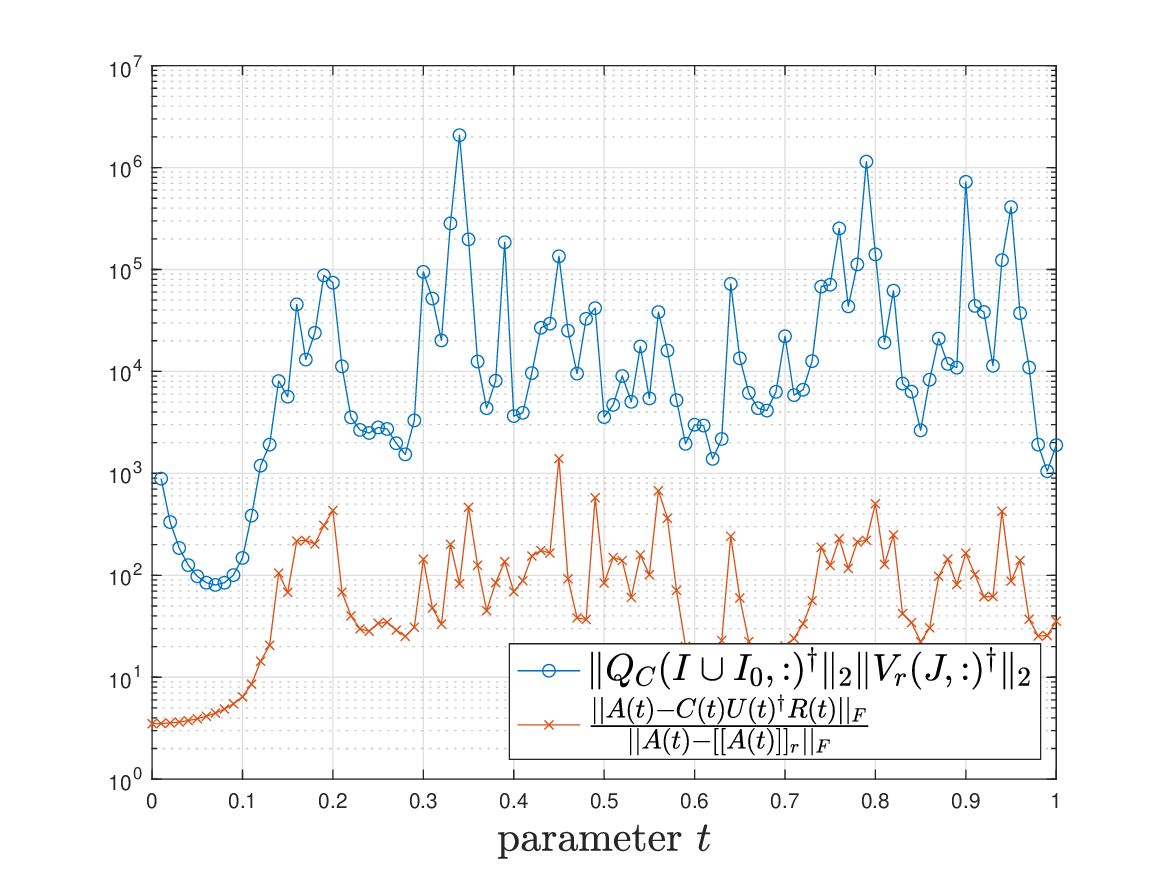}} 
\centering 
\caption{Testing the gap between theory and practice for the CUR bound in \eqref{eq:CAboundOS} using the same set of indices for every parameter value. The parameter-dependent matrix in this example is the same as the one in Section \ref{subsec:RoleTol} (Equation \eqref{eq:synprob}). The target rank is $r = 33$ in this experiment.}
\label{fig:SamePlot}
\end{figure}

The results are shown in Figure~\ref{fig:SamePlot}. We begin with an initial set of indices and keep using the same set of indices for all other parameter values. The target rank is $33$ in this experiment. As shown in Figure~\ref{subfig:SameP}, despite using the same set of indices, we achieve a relative error of about $10^{-7}$ throughout, losing only about $1$--$2$ digits of accuracy compared to the initial accuracy. On the other hand, in Figure~\ref{subfig:SameC}, we observe that the bound provided by \eqref{eq:CAboundOS} significantly overestimates the true bound. This demonstrates that the quantity $\norm{Q_C(I \cup I_0,:)^\dagger}_2 \norm{V_r(J,:)^{-1}}_2$ alone is insufficient to explain the effectiveness of the CUR decomposition. Such gaps between theoretical bounds and practical performance are common and have been observed in other works, such as \cite{SorensenEmbree2016}.

In many problems, $\delta$ may be large, yet the CUR decomposition can provide a far better approximation than what is suggested by the bound in \eqref{eq:CAboundOS}. 
Moreover, in most practical cases $\norm{A-\lowrank{A}{r}}_F$ is unknown. 
This observation motivates us to check whether the previously obtained indices can still be used before recalculating them entirely. To guard against potential large errors, we incorporate error and rank estimation in our algorithms, which we discuss in detail in Sections \ref{subsec:mainalg} and \ref{subsec:cheapalg}.

The second observation is in the role that $I_0$ plays in \eqref{eq:CAboundOS}. The set of indices $I_0$ oversamples the row indices to improve the term $\norm{Q_C(I\cup I_0,:)^\dagger}_2$, benefiting from the fact that $\norm{Q_C(I\cup I_0,:)^\dagger}_2 \leq \norm{Q_C(I,:)^{-1}}_2$, which follows from the Courant-Fischer min-max theorem. Furthermore, when $|I_*|>|J|$, the core matrix $A(I_*,J)$ has larger singular values than when $A(I,J)$ is a square matrix. This improves the accuracy and stability in the computation of the CUR decomposition; see \cite{pn24cur} for a detailed discussion. The concept of oversampling has been explored in prior works such as \cite{AndersonDuMahoneyMelgaardWuGu2015,DonelloPalkarNaderiDelReyFernandezBabaee2023,pn24cur,PeherstorferDrmacGugercin2020,ZimmermannWillcox2016} and it is known that oversampling improves the accuracy of the CUR decomposition. In light of this observation, we adopt oversampling, and for definiteness choose to oversample rows in this work, that is, $|I_0|>0$. We summarize a version of the oversampling algorithm from \cite{pn24cur} in Algorithm \ref{alg:oversample}, which increases the minimum singular value(s) of $Q_C(I,:)$ by finding unchosen indices that enrich the trailing singular subspace of $Q_C(I,:)$.

\begin{algorithm}[!ht]\footnotesize
  \caption{Oversampling for CUR}
  \label{alg:oversample}
  \begin{algorithmic}[1]
\Require{$A\in \R^{m\times n}$, column indices $J$, row indices $I$, with $|I|=|J|= r< m,n$, oversampling parameter $p (\leq r)$}
\Ensure{Extra indices $I_0$ with $|I_0| = p$}
\vspace{0.5pc}
\Statex $I_0 = \mathtt{OS}(A,I,J,p)$
\State $[Q_C,\sim] = \mathtt{qr}\left(A(:,J)\right)$,
\State $[\sim,\sim,V] = \mathtt{svd}(Q_C(I,:))$,
\State Set $V_{-p} = V(:,r-p+1:r)$, the trailing $p$ right singular vectors of $Q_C(I,:)$.
\State Set $Q_{-C} = Q_C([m]-I,:)V_{-p}$,
\State Apply CPQR on $Q_{-C}^T$. Let $I_0$ be the $p$ extra row pivots.
\end{algorithmic}
\end{algorithm}

Algorithm \ref{alg:oversample} obtains extra row indices $I_0$, distinct from $I$, by trying to increase the minimum singular value of $Q_C(I,:)$. The algorithm does this by first projecting the trailing singular subspace of $Q_C(I,:)$ onto $Q_C$ and choosing, through pivoting, the unchosen indices that contribute the most to the desired subspace, thereby increasing the minimum singular value of $Q_C(I,:)$. The details along with its motivation using cosine-sine decomposition can be found in \cite{pn24cur}. The complexity of Algorithm \ref{alg:oversample} is $\bigO(mr^2+nrp)$ where the dominant costs come from computing the QR decomposition $\bigO(mr^2)$ in line $1$ and the cost of matrix-matrix multiply $\bigO(nrp)$ in line $4$.

\subsection{Computing indices from scratch} \label{subsec:scratch} 
In the next two sections, we discuss our algorithms, \acur{} and \qcur{}. For both algorithms, we need to start with an index set for the CUR decomposition. This can sometimes be done offline, for example, if we are solving a Matrix PDE, we can compute the initial set of indices for the CUR decomposition from the initial condition matrix. In this work, we use Algorithm \ref{alg:pivotrankest} to get a set of indices from scratch. The procedure is to first approximate the $(\epsilon/\sqrt{n})$-rank of the initial matrix $A(t_1)$ using randomized rank estimation and then apply pivoting on a random sketch with the estimated $(\epsilon/\sqrt{n})$-rank to get the initial set of indices. Here, we approximate the $(\epsilon/\sqrt{n})$-rank to ensure a relative error of $\epsilon$ for the CUR decomposition in the Frobenius norm. More specifically, the $1/\sqrt{n}$ factor arises from the worst case scenario where the trailing singular values of $A(t)$ are all equal. In such cases, the $(\epsilon/\sqrt{n})$-rank is required to guarantee a relative accuracy of $\epsilon$ in the Frobenius norm. If additional information about the spectrum of $A(t)$ is available--for instance, if $A(t)$ exhibits rapid singular value decay--then $\epsilon/\sqrt{n}$ can be replaced by $C\epsilon$, where $C<1$ is a constant. However, to guarantee accuracy, we use $(\epsilon/\sqrt{n})$-rank throughout work. The procedure for computing indices from scratch will be used for getting an initial set of indices for \acur{} and \qcur{}, but also in \acur{} when the error becomes too large that we need to recompute the indices altogether from scratch.

\subsection{\acur{} algorithm for accuracy} \label{subsec:mainalg}
\acur{} is an algorithm aimed at efficiently computing an accurate low-rank approximation of a parameter-dependent matrix with a given error tolerance. An overview of \acur{} goes as follows. The algorithm starts by computing the initial set of indices for $A(t_1)$ as discussed in Section \ref{subsec:scratch}. For subsequent parameter values, i.e., $A(t_j)$ for $j = 2,...,q$, we first verify whether the indices obtained from the previous parameter value are able to meet some given error tolerance, using an estimate of $\|A-CU^\dagger R\|_F$. If so, we continue to the next parameter value. Should the indices fail to meet the tolerance, we make low-cost minor modifications to the indices and test whether the adjustments satisfy the error tolerance. If this is still unsuccessful, we recompute the indices entirely from scratch and move onto the next parameter value. Minor modifications include adding a small fixed number of indices, reordering them based on their importance and removing some if necessary. The details are in the following paragraph. The addition and removal of indices, as well as recomputation of indices using rank estimation if necessary, make the algorithm rank-adaptive and computing the relative error at each instance make the algorithm certifiable. \acur{} is presented in Algorithm \ref{alg:main}.

\begin{flushleft}
\begin{minipage}[!h]{\textwidth}
\renewcommand*\footnoterule{}
\begin{savenotes}
\begin{algorithm}[H]\footnotesize
  \caption{\acur{} with certified accuracy}
  \label{alg:main}
  \begin{algorithmic}[1]
\Require{Parameter-dependent matrix $A(t)\in \R^{m\times n}$ ($m\geq n$), evaluation points $t_1,...,t_q$, error tolerance $\epsilon$, error sample size $s$ (rec. $s=5$), oversampling parameter $p$.}
\Ensure{Matrices $C(t_j), R(t_j), U(t_j)$ defining a subset of columns and rows of $A(t_j)$ and their intersection for $j = 1,...,q$.}
\vspace{0.5pc}
\State $[I,J] = \mathtt{Rand\_ Pivot \_ RankEst}(A(t_1),\epsilon/\sqrt{n})$ \Comment{Algorithm \ref{alg:pivotrankest}}
\State $I_0 = \mathtt{OS}(A(t_1),I,J,p)$ \Comment{Algorithm \ref{alg:oversample}}
\State Set $C(t_1) = A(t_1)(:,J),$ $R(t_1) = A(t_1)(I\cup I_0,:)$ and $U(t_1) = A(t_1)(I\cup I_0,J)$,
\For{$j = 2,...,q$}
\State Draw a Gaussian embedding $\Gamma_s \in \R^{s\times m}$
\State Sketch $X_s = \Gamma_s A(t_j) \in \R^{s\times n}$
\State Set $E_s = X_s-\Gamma_s A(t_j)(:,J)A(t_j)(I\cup I_0,J)^\dagger A(t_j)(I\cup I_0,:)$ \Comment{Sketch of the residual}
\State Approximate relative error $E = \norm{E_s}_F/\norm{X_s}_F$
\If{$E>\epsilon$}
\State $[I_1,J_1] = \mathtt{Pivot\_Residual}(E_s,I,J)$ \Comment{Algorithm \ref{alg:pivotRes}}
\State Append indices $I \leftarrow I \cup I_1$ and $J \leftarrow J\cup J_1$
\State $[\sim,\sim,I_2] = \mathtt{sRRQR}\left(A(t_j)\left(I\cup I_0,J\right)^T\right)$ \footnote{$[Q,R,I] = \mathtt{sRRQR}(A)$ performs strong rank-revealing QR factorization \cite{GuEisenstat1996} on $A$ where $Q$ is an orthonormal matrix, $R$ is an upper triangular matrix and $I$ is a set of pivots.}
\State $[\sim,R,J_2] = \mathtt{sRRQR}\left(A(t_j)\left(I\cup I_0,J\right)\right)$
\State $r = \left\vert\left\{R_{ii}: |R_{ii}| > \epsilon/\sqrt{n} \cdot |R_{11}|\right\}\right\vert$ \Comment{Estimate relative $(\epsilon/\sqrt{n})$-rank}
\State $I \leftarrow (I\cup I_0)\left(I_2(1:r)\right)$ \Comment{Select $r$ important row indices from $I\cup I_0$}
\State $I_0 \leftarrow (I\cup I_0)\left(I_2(r+1:r+p)\right)$ \Comment{Select the next $p$ important row indices from $I\cup I_0$}
\State $J \leftarrow J\left(J_2(1:r)\right)$ \Comment{Select $r$ important column indices from $J$}
\State Set $E_s = X_s-\Gamma_s A(t_j)(:,J)A(t_j)(I\cup I_0,J)^\dagger A(t_j)(I \cup I_0,:)$
\State Recompute relative error $E = \norm{E_s}_F/\norm{X_s}_F$
\If{$E>\epsilon$} 
\State $[I,J] = \mathtt{Rand\_ Pivot \_ RankEst}(A(t_j),\epsilon/\sqrt{n})$ \Comment{Recompute indices (Algorithm \ref{alg:pivotrankest})} 
\State $I_0 = \mathtt{OS}(A(t_j),I,J,p)$ \Comment{Algorithm \ref{alg:oversample}}
\EndIf
\EndIf
\State Set $C(t_j) = A(t_j)(:,J),$ $R(t_j) = A(t_j)(I\cup I_0,:)$ and $U(t_j) = A(t_j)(I\cup I_0,J)$.
\EndFor
\end{algorithmic}
\end{algorithm}
\end{savenotes}
\end{minipage}
\end{flushleft}
\vspace{0.2cm}

\acur{} (Algorithm \ref{alg:main}) takes the parameter-dependent matrix $A(t)$, evaluation points $t_1,...,t_q \in D$, error tolerance $\epsilon$, error sample size $s$, and oversampling parameter $p$ as input. The algorithm produces the low-rank factors $C(t_j), R(t_j), U(t_j)$ as output such that the factors are a subset of columns, rows and their intersection of $A(t_j)$, respectively, and $A(t_j) \approx C(t_j)U(t_j)^\dagger R(t_j)$. If the entries of the parameter-dependent matrix can be computed quickly afterwards, the algorithm can be modified to output only row and column indices, saving storage. \acur{} is quite involved, so we break it down into two core parts: 
\begin{itemize}
    \item Section \ref{step:ii}: Error estimation using previous indices (lines $5$--$8$),
    \item Section \ref{step:iii}: Low-cost minor modifications (lines $10$--$19$).
\end{itemize}
The other lines involve computing the indices from scratch, which follow the same discussion as Section \ref{subsec:scratch}.

\subsubsection{ Initial error estimation} \label{step:ii}
In the for loop, the first task is to quantify how well the previous set of indices perform on the current parameter value. Error estimation is used for this task where we generate a Gaussian embedding $\Gamma_s$ in line $5$, sketch the matrix corresponding to the current parameter value $X_s = \Gamma_s A(t_j)$ in line $6$ and the residual error matrix for the CUR decomposition in line $7$, and finally, estimating the relative error $E$ using the sketches in line $8$. The sketch of the residual error matrix $E_s$ is used for two tasks: approximating the relative error (lines $8,19$) and computing additional indices to reduce the relative error (line $10$). If the relative error $E$ is less than the tolerance $\epsilon$, we use the previous set of indices and store the corresponding columns, rows and their intersection in line $23$, after which we continue to the next parameter value $t_j$. On the other hand, if the relative error exceeds the tolerance, we make low-cost minor modifications in an attempt to reduce the relative error, as we describe next.

\subsubsection{ Low-cost minor modifications} \label{step:iii}
The low-cost modifications involve trying to enlarge the index set using the sketch of the residual error matrix $E_s$ (lines $10$--$11$), reordering them by importance (lines $12$--$13$), removing unimportant indices (lines $14$--$17$), and recomputing the relative error (lines $18$--$19$). More specifically, in line $10$, randomized pivoting is used on the sketch of the residual matrix, $E_s$ as in Algorithm \ref{alg:pivotRes} to obtain an additional set of $s$ row and column indices. We append the extra set of $s$ indices to the original sets of indices $I$ and $J$ in line $11$, and order the indices in terms of their importance in lines $12$--$13$. Here, Gu-Eisenstat's strong rank-revealing QR factorization \cite{GuEisenstat1996} is used, but we can use other strong rank-revealing algorithms or more expensive algorithms with strong theoretical guarantees such as \cite{CortinovisKressner2020,Osinsky2023}, because we are working with a small $O(r)\times r$ matrix $A(t_j)(I\cup I_0,J)$. In line $14$, we approximate the $(\epsilon/\sqrt{n})$-rank of $A(t_j)(I\cup I_0,J)$ by looking at the diagonal entries of the upper triangular factor in sRRQR as they give an excellent approximation to the singular values of the original matrix.\footnote{Instead of approximating the $(\epsilon/\sqrt{n})$-rank of $A(t_j)(I\cup I_0,J)$, we can set the rank tolerance slightly smaller, for example, $0.5\epsilon/\sqrt{n}$, to account for the approximation error and decrease the chance of recomputing the indices from scratch (lines $21$--$22$) by keeping slightly more indices than necessary.} In lines $15$--$17$, we truncate based on the rank computed in line $14$. The ordering of the indices is important here. We recompute the relative error in lines $18$--$19$,\footnote{For theoretical guarantee, the Gaussian embedding $\Gamma_s$ needs to be independent from the modified indices \cite{GrattonTitley-Peloquin2018}, however the extra set of indices is dependent on $\Gamma_s$ as we apply pivoting on the sketched residual in line $10$. Nonetheless, the dependency is rather weak so $\Gamma_s$ can be reused without losing too much accuracy.} and if the relative error $E$ is smaller than the tolerance, we store the new set of columns, rows of $A(t_j)$ and their intersection in line $23$ and continue to the next parameter value. If the relative error exceeds the tolerance, there are two recourses. \acur{} outlines the first option in lines $21$--$22$ where we recompute the indices altogether from scratch. In line $21$, we perform both randomized rank estimation and randomized pivoting for two reasons. First, randomized rank estimation is extremely cheap when we have the row sketch already and second, when there is a sudden large change in rank, we can use randomized rank estimation to adjust the rank efficiently; note that when minor modifications did not help in lines $10$--$17$, the matrix may have changed drastically.

\paragraph{An alternative to recomputing the indices}
 An alternative approach, not included in \acur{}, is to increase the error sample size $s$ gradually by incrementing its value, say $2s,4s,$ and so on, and updating the sketch accordingly by appending to it. As $s$ increases, the number of additional important indices obtained from randomized pivoting in line $10$ also increases, which will eventually reduce the relative error to less than the tolerance $\epsilon$.

\paragraph{Role of the error sample size $s$}
The value of $s$ can be important for two reasons: error estimation and minor modifications. If we set $s$ to be larger, we obtain higher accuracy in our error estimation at the cost of more computation. However, we obtain a larger set of extra indices for the minor modification. This can help the algorithm avoid the if statement in line $20$, which recomputes the indices from scratch; see also Section \ref{subsec:RoleS}. We can also create a version with increasing $s$ if minor modifications fail to decrease the relative error below the tolerance, as described in the final part of the previous paragraph.

\paragraph{Complexity} Let $r_j$ be the size of the index set for parameter $t_j$, i.e., the rank of the CUR decomposition for $A(t_j)$, $Q_1$ the set of parameter values for which we recompute the indices from scratch, i.e., invoked lines $1,2,21,22$, and $Q_2$ the set of parameter values for which we do not need to recompute the indices from scratch. Note that $t_1 \in Q_1$, $Q_1 \cup Q_2 = \{t_1,...,t_q\}$ and $Q_1 \cap Q_2 = \emptyset$.  Then the complexity of \acur{} is 
\begin{equation}
    \bigO\Big(\sum\limits_{j \in Q_1} \left(r_j T_{A(t_j)} + (m+n)r_j^2+nr_j (r_j+p)\right) + \sum\limits_{j \in Q_2} \left(s T_{A(t_j)} + (m+n)r_j s+ n s p\right)\Big)
\end{equation} where $T_{A(t_j)}$ is the cost of matrix-vector multiply with $A(t_j)$ or $A(t_j)^T$. The dominant cost, which can be quadratic with respect to $m$ and $n$, comes from sketching, which cost $\bigO(r_j T_{A(t_j)})$ in lines $1,21$ and $\bigO(sT_{A(t_j)})$ in line $6$. The dominant linear costs come from pivoting and oversampling in lines $1,2,21,22$, which cost $\bigO((m+n)r_j^2+nr_jp)$ and computing the sketch of the residual error matrix in lines $7,18$, which is $\bigO(mr_j s + ns(r_j+p))$. All the other costs are smaller linear costs or sublinear with respect to $m$ and $n$, for example, the SRTT cost in randomized rank estimation in lines $1,21$ is $\bigO(nr_j\log r_j)$ and the cost of sRRQR in lines $12,13$ is $\bigO((r_j+p+s)(r_j+s)^2)$. When $s,p = \bigO(1)$ the complexity simplifies to
\begin{equation} \label{eqn:simpcomp}
    \bigO\left(\sum\limits_{j \in Q_1} \left(r_j T_{A(t_j)} + (m+n)r_j^2\right) + \sum\limits_{j \in Q_2} \left(T_{A(t_j)} + (m+n)r_j\right)\right).
\end{equation} 

It turns out that in many applications, the cardinality of the set $Q_1$ is small when we choose the oversampling parameter $p$ and the error sample size $s$ appropriately, making the algorithm closer to $\bigO(T_{A(t_j)})$ for the $j$th parameter value; see Sections \ref{subsec:RoleP} and \ref{subsec:RoleS}. This makes \acur{} usually faster than other existing methods such as dynamical low-rank approximation \cite{LubichOseledets2014}, and randomized SVD and generalized Nystr\"om method \cite{KressnerLam2023}; see Section \ref{subsec:Speed}.

\subsection{\qcur{} algorithm for speed} \label{subsec:cheapalg}
The bottleneck in \acur{} is usually in the sketching, which is crucial for reliably computing the set of row and column indices for the CUR decomposition and for error estimation, ensuring the algorithm's robustness. Without sketching, our algorithm would have linear complexity with respect to $m$ and $n$; see \eqref{eqn:simpcomp}. In this section, we propose \qcur{} that achieves linear complexity after the initial phase of computing the initial set of indices for $A(t_1)$. \qcur{} offers the advantage that the number of rows and columns required to be seen for each parameter value is approximately its target rank. Therefore, there is no need to read the whole matrix, which is advantageous in settings where the entries are expensive to compute. Unfortunately, the fast version does suffer from adversarial examples as \qcur{} does not look at the entire matrix; see Section \ref{subsec:adversarial}. Nevertheless, the algorithm seems to perform well in practice for easier problems as illustrated in Section \ref{sec:numill}.

An overview of \qcur{} goes as follows. The algorithm starts by computing the initial set of indices for $A(t_1)$ as discussed in Section \ref{subsec:scratch}. Then for each subsequent parameter value, we first form the core matrix $U_j = A(t_j)(I,J)$ where $I$ and $J$ include extra indices from the buffer space and oversampling from the previous step $j-1$. Here, the buffer space is the extra indices kept to quickly detect a possible change in rank or the important indices as we iterate through the parameter values. We then order the indices by importance and compute the ($\epsilon/\sqrt{n}$)-rank of the core matrix $U_j$.\footnote{We compute the ($\epsilon/\sqrt{n}$)-rank of the core matrix to aim for a relative low-rank approximation error of $\epsilon$ in the Frobenius norm. As in \acur{}, we can set the rank tolerance slightly smaller, say $0.5\epsilon/\sqrt{n}$ to account for the approximation error.} The order of indices is important here as we explain in (ii) below. If the rank increased from the previous iteration, we add extra indices from the buffer space and replenish it by invoking the oversampling algorithm using $\mathtt{OS}(A,I,J,\delta_r)$ (Algorithm \ref{alg:oversample}) where $\delta_r$ is the rank increase. Conversely, if the rank decreases, we simply remove some indices. Finally, we store the corresponding rows and columns of $A(t_j)$ and their intersection and move to the next parameter value. \qcur{} is presented in Algorithm \ref{alg:cheap}.

\begin{algorithm}[h]\footnotesize
  \caption{\qcur{} for speed}
  \label{alg:cheap}
  \begin{algorithmic}[1]
\Require{Parameter-dependent matrix $A(t)\in \R^{m\times n}$ ($m\geq n$), evaluation points $t_1,...,t_q$, buffer space $b$, oversampling parameter $p$, rank tolerance $\epsilon$.}
\Ensure{Matrices $C(t_j), R(t_j), U(t_j)$ defining a subset of columns and rows of $A(t_j)$ and their intersection for $j = 1,...,q$.}
\vspace{0.5pc}
\State $[I,J] = \mathtt{Rand\_ Pivot\_RankEst}(A(t_1),\epsilon/\sqrt{n})$ \Comment{Algorithm \ref{alg:pivotrankest}}
\State Set $r\leftarrow |I|$
\State $I_0 = \mathtt{OS}(A(t_1),I,J,p+b)$ \Comment{Algorithm \ref{alg:oversample}; includes buffer space and oversampling}
\State $J_0 = \mathtt{OS}(A(t_1)^T,J,I,b)$ \Comment{Algorithm \ref{alg:oversample}; includes buffer space}
\State Set $I \leftarrow I\cup I_0$ and $J \leftarrow J\cup J_0$
\State Set $C(t_1) = A(t_1)(:,J(1:r)), R(t_1) = A(t_1)(I(1:r+p),:), U(t_1) = A(t_1)(I(1:r+p),J(1:r))$
\For{$j = 2,...,q$}
\State Set $U = A(t_j)(I,J) \in \R^{(r+b+p)\times (r+b)}$
\State $[\sim,\sim,I_1] = \mathtt{sRRQR}\left(U^T\right)$ 
\State $[\sim,R,J_1] = \mathtt{sRRQR}\left(U\right)$
\State $r_0 = \left\vert\left\{R_{ii}: |R_{ii}| > \epsilon/\sqrt{n}\cdot |R_{11}|\right\}\right\vert$ \Comment{Estimate relative  $\left(\epsilon/\sqrt{n}\right)$-rank}
\If{$r_0 \leq r$}
\State Set $I \leftarrow I\left(I_1(1:r_0+b+p)\right)$ \Comment{Reorder row indices by importance and truncate}
\State Set $J \leftarrow J\left(J_1(1:r_0+b)\right)$ \Comment{Reorder column indices by importance and truncate}
\Else 
\State $I_2 = \mathtt{OS}(A(t_j),I,J,r_0-r)$ \Comment{Replenish row indices}
\State $J_2 = \mathtt{OS}(A(t_j)^T,J,I,r_0-r)$ \Comment{Replenish column indices}
\State Set $I \leftarrow I(I_1) \cup I_2$ and $J \leftarrow J(J_1) \cup J_2$ \Comment{Reorder and add indices}
\EndIf
\State Set $r\leftarrow r_0$ \Comment{Update $\left(\epsilon/\sqrt{n}\right)$-rank}
\State Set $C(t_j) = A(t_j)(:,J(1:r)), R(t_j) = A(t_j)(I(1:r+p),:), U(t_j) = A(t_j)(I(1:r+p),J(1:r))$
\EndFor
\end{algorithmic}
\end{algorithm}


The main distinction between \acur{} and \qcur{} lies in the existence of error estimation. In \acur{}, an error estimate is computed for each parameter value, which ensures that the row and the column indices are guaranteed to be good with high probability. This is absent in \qcur{}. Consequently, if there is a sudden change in rank or the previous set of indices are unimportant for the current parameter value, \qcur{} may deliver poor results. To mitigate this disadvantage, \qcur{} incorporates a buffer space (i.e., additional indices in $I,J$) as input to assist with detecting changes in rank. \qcur{} involves many heuristics, so we break it down to discuss and justify each part of the algorithm. We divide \qcur{} into three parts: 
\begin{itemize}
    \item Section \ref{step:fast1}: Rank estimation using the core matrix (lines $8$--$11$),
    \item Section \ref{step:fast2}: Rank decrease via truncation (lines $12$--$14$),
    \item Section \ref{step:fast3}: Rank increase via oversampling (lines $15$--$18$). 
\end{itemize}
The initial set of indices are constructed from scratch as in Section \ref{subsec:scratch}, with the indices for the buffer space and oversampling obtained using the oversampling algorithm from Algorithm \ref{alg:oversample}.

\subsubsection{Rank estimation using the core matrix} \label{step:fast1}
Upon entering the for loop in line $7$, we first compute the core matrix $U = A(t_j)(I,J)\in \R^{(r+b+p)\times (r+b)}$ (line $8$) for the current parameter value, which includes extra indices from the buffer space and oversampling. As \qcur{} prioritizes efficiency over accuracy, the core matrix is a sensible surrogate to approximate the singular values of $A(t_j)$ without looking at the entire matrix. The extra indices from the buffer space and oversampling will assist with the approximation as the buffer space allows the core matrix to approximate more singular values, helping to detect a potential $(\epsilon/\sqrt{n})$-rank change. The additional set of row indices from oversampling improves the accuracy of the estimated singular values of $A(t_j)$ (line $10$) by the Courant-Fischer min-max theorem, which is used for approximating the relative $(\epsilon/\sqrt{n})$-rank in line $11$. In addition, oversampling increases the singular values of the core matrix, thereby allowing \qcur{} to increase the rank when appropriate and reduce the error. However, we see in the experiments (Section \ref{subsec:RoleP}) that the role that oversampling plays in \qcur{} is rather complicated. We use sRRQR to order the columns and rows of the core matrix by importance in lines $9$ and $10$, and subsequently estimate the $(\epsilon/\sqrt{n})$-rank of the core matrix $U$ using the diagonal entries of the upper triangular factor in sRRQR in line $11$. We use the $(\epsilon/\sqrt{n})$-rank to aim for a tolerance of $\epsilon$ in the low-rank approximation. If the computed rank, $r_0$ is smaller than the previous rank, $r$, then we decrease the index set via truncation. Otherwise, we increase the index set via oversampling. The details are discussed below.

\subsubsection{ Rank decrease via truncation}\label{step:fast2}
After estimating the $(\epsilon/\sqrt{n})$-rank $r_0$, if $r_0$ has not increased from the previous iteration, we adjust the number of indices (either decrease or stay the same) in lines $13$ and $14$ by applying the order of importance computed using sRRQR in lines $9$ and $10$, and truncating, if necessary. Specifically, we truncate the trailing $r-r_0$ row and column indices. At the end of this procedure, the algorithm ensures $|I| = r_0+b+p$ and $|J| = r_0+b$.

\subsubsection{Rank increase via oversampling} \label{step:fast3}
If the rank has increased, i.e. $r_0 > r$, we refill the extra indices by the amount the rank has increased by, i.e., $r_0-r$. Unlike in \acur{}, for which we have the sketched residual matrix, we do not have a good proxy to obtain a good set of extra indices. Therefore, we use the already-selected columns and rows to obtain an extra set of column and row indices.\footnote{The only recourse to obtaining a reliable set of extra indices is to view the entire matrix; see \cite[\S~17.5]{martinsson2019fast} for a related discussion. However, this puts us in a similar setting as \acur{}, which we recommend when ensuring accuracy is the priority rather than speed.} We use the oversampling algorithm (Algorithm \ref{alg:oversample}) to achieve this as it only requires the selected rows and columns as input. Adding extra indices using oversampling has been suggested before in \cite{DonelloPalkarNaderiDelReyFernandezBabaee2023} to make small rank increases in the context of dynamical low-rank approximation. We get the extra indices from oversampling in lines $16$ and $17$ and append to the original set of indices in line $18$. At the end of this procedure, the algorithm ensures $|I| = r_0+b+p$ and $|J| = r_0+b$.

\qcur{} relies on heuristics to make sensible decisions for estimating the $(\epsilon/\sqrt{n})$-rank, selecting the important indices, and adding indices as we describe above. Despite \qcur{} being susceptible to adversarial examples, as we see later in Section \ref{subsec:adversarial}, it should work on easier problems. For example, if the singular vectors of $A(t)$ are incoherent \cite[\S~9.6]{MartinssonTropp2020} or (approximately) Haar-distributed, \qcur{} will perform well because uniformly sampled columns and rows approximate $A(t)$ well \cite{ChiuDemanet2013}. In such cases, rank adaption is also expected to work well, as the submatrix of $A(t)$ with oversampling behaves similarly to the two-sided sketch used for randomized rank estimation \cite{MeierNakatsukasa2024}; see Section \ref{subsec:rankest}. Additionally, with the aid of buffer space, \qcur{} is expected to detect rank changes occurring in the problem and efficiently adapt to the rank changes using oversampling or truncation. Therefore, for problems where $A(t)$ is incoherent for the parameter values of our interest, \qcur{} is expected to perform effectively.

\paragraph{Role of the buffer space $b$}
The main role of the buffer space $b$ is to identify rank changes while iterating through the parameter values. A larger value of $b$ enhances the algorithm's ability to accurately detect rank changes at an expense of increased complexity. Furthermore, a larger buffer size allows the algorithm to make larger rank changes as the maximum possible rank change at any point in the algorithm is $b$; see Section \ref{subsec:RoleS}.

\paragraph{Complexity} Let $r_j$ be the size of the index set for parameter $t_j$, i.e., the rank of the CUR decomposition for $A(t_j)$. The dominant cost $\bigO\left(r_1 T_{A(t_1)}\right)$ in \qcur{} comes from line $1$ for computing the initial set of indices for $A(t_1)$ excluding the buffer space and oversampling. Besides the computation of initial indices, the complexity is at most linear in $m$ and $n$. For $j = 2,...,q$, the complexity of \qcur{} for the $j$th parameter value is $\bigO((m+n)(r_j+p+b)^2)$ if the rank has increased from $t_{j-1}$ to $t_j$. The cost involves the usage of oversampling to replenish extra indices. If the rank has either decreased or stayed the same then the complexity is $\bigO((r_j+p+b)(r_j+b)^2)$ for using sRRQR. Hence, excluding the first line, \qcur{} runs at most linearly in $m$ and $n$ for each parameter value, making it remarkably fast as demonstrated in Section \ref{subsec:Speed}.

\section{Numerical Illustration} \label{sec:numill}
In this section, we illustrate \acur{} (Algorithm \ref{alg:main}) and \qcur{} (Algorithm \ref{alg:cheap}) using numerical experiments. The code for the strong rank-revealing QR algorithm \cite{GuEisenstat1996} is taken from \cite{sRRQRcode}. All experiments were performed in MATLAB version 2021a using double precision arithmetic on a workstation with Intel® Xeon® Silver 4314 CPU @ 2.40GHz ($16$ cores) and 512GB memory.

In the following sections, we perform different numerical simulations with varying input parameters to test \acur{} and \qcur{}. These simulations allow us to test the following.
\begin{enumerate}
    \item The accuracy of the approximation obtained from the two algorithms by varying the tolerance $\epsilon$,
    \item The role that the oversampling parameter $p$ plays in our algorithms,
    \item The role that the error sample size $s$ plays in \acur{} and the buffer size $b$ in \qcur{},
    \item How our algorithms perform on adversarial problems,
    \item The speed of our algorithms against existing methods in \cite{KressnerLam2023,LubichOseledets2014}.
\end{enumerate} It is worth pointing out that in the first three examples $\min_j \norm{A(t_{j})-A(t_{j-1})}_2 > 0.03$ for all $j$, indicating that the performance of \acur{} and \qcur{} is not explained by the fact that the matrices are undergoing small perturbations as we iterate through the parameters.

In the experiments, we measure how many times \acur{} needs to enter some of the `if' statements, as the heavier computations are usually done inside them. This happens when the current set of indices do not meet the given tolerance in \acur{}. We define the following
\begin{itemize}
    \item $h_1$ is the number of times \acur{} has invoked lines $10$--$18$, but not lines $20$--$21$. This happens when only minor modifications were needed to meet the tolerance. The complexity is $\bigO(T_{A(t_j)}+(m+n)r_j^2)$ for the $j$th parameter value.
    \item $h_2$ is the number of times \acur{} has invoked lines $20$--$21$. This happens when minor modifications were not enough to meet the tolerance and the row and the column indices had to be recomputed from scratch. The complexity is $\bigO(r_j T_{A(t_j)}+(m+n)r_j^2)$ for the $j$th parameter value.
\end{itemize} In the plots that depict the $(\epsilon/\sqrt{n})$-rank of the low-rank approximation, a large dot is used to indicate the parameter value for which the indices were recomputed from scratch in \acur{}. The number of large dots is equal to $h_2$. In \qcur{}, besides the computation of initial sets of indices, the heaviest computation is done when the rank has increased; i.e. in lines $16$ and $17$ when oversampling is invoked. However, this is heavily dependent on the given problem, so we do not track this. For reference, in \qcur{}, a large dot is used to indicate the parameter value where oversampling was applied when the rank increased.

\subsection{Varying the tolerance $\epsilon$} \label{subsec:RoleTol}
In this section, we test the accuracy of our algorithms by varying the tolerance parameter $\epsilon$. We use the same synthetic example \cite{CerutiLubich2022} used for Figure \ref{fig:SamePlot}, which is given by
    \begin{equation} \tag{\ref{eq:synprob}}
    A(t) = e^{tW_1}e^t D e^{tW_2}, t\in [0,1].
\end{equation} This example tests robustness against small singular values.

\begin{figure}[!ht]
\subfloat[\acur{} with $p = 5$, $s = 5$]{\label{subfig:TolM1}\includegraphics[scale = 0.33]{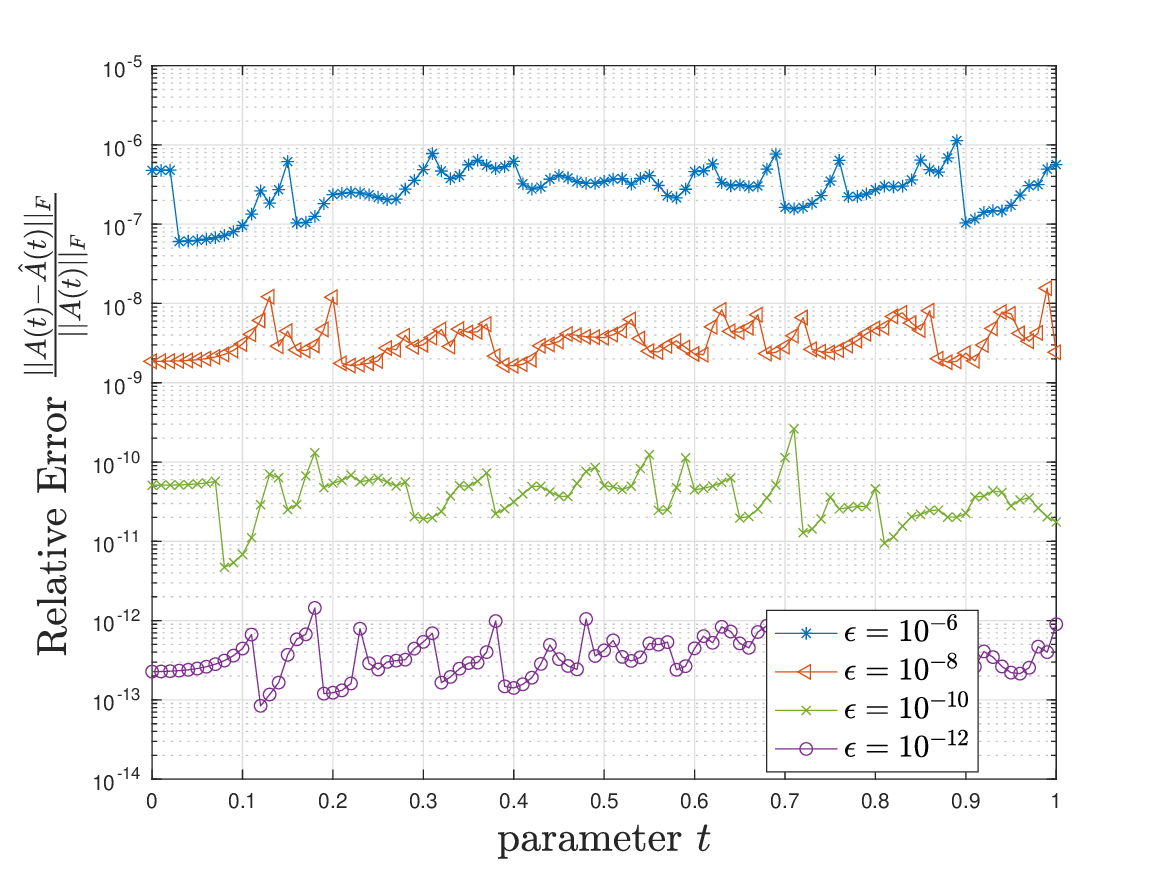}} \hspace{-0.5cm}
\subfloat[Rank changes for \acur{}]{\label{subfig:TolM2}\includegraphics[scale = 0.33]{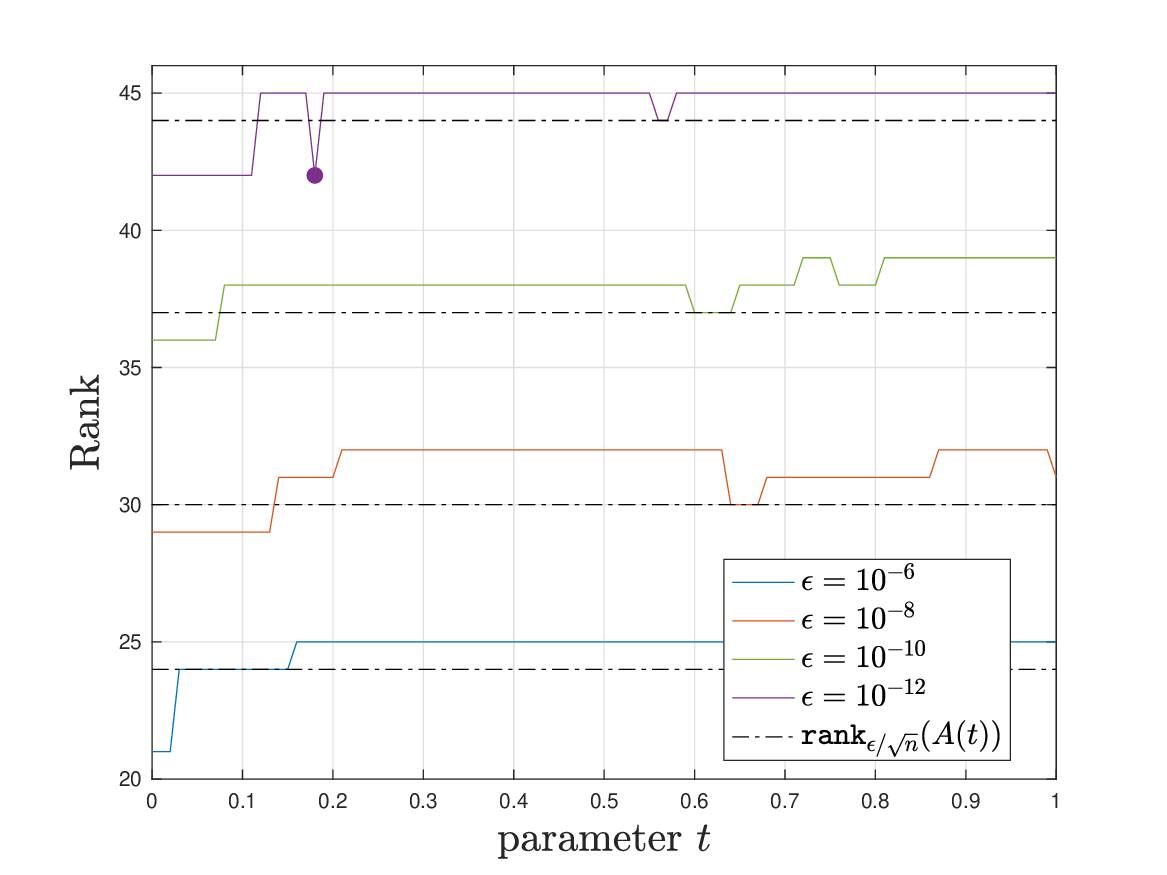}} \vspace{-0.4cm} \\ 
\subfloat[\qcur{} with $p= 5$, $b = 5$]{\label{subfig:TolA1}\includegraphics[scale = 0.33]{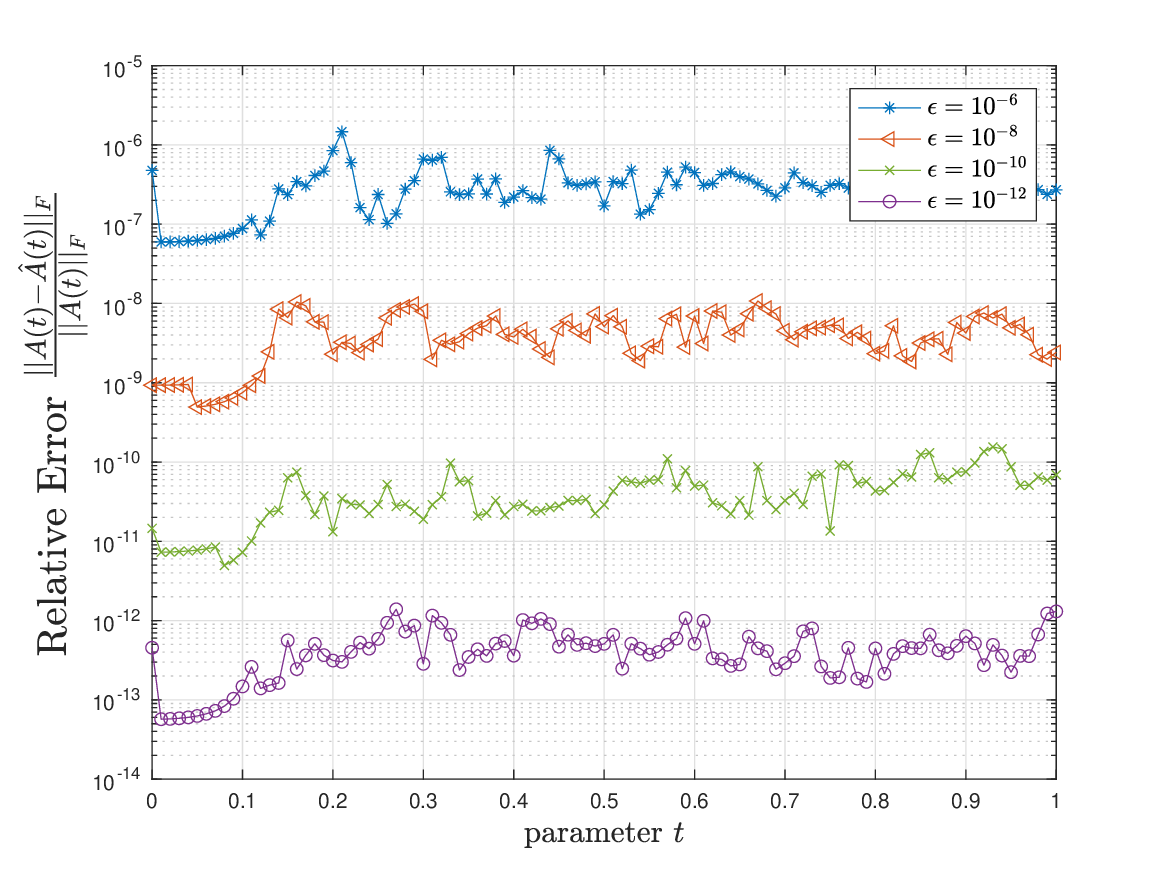}} \hspace{-0.5cm}
\subfloat[Rank changes for \qcur{}]{\label{subfig:TolA2}\includegraphics[scale = 0.33]{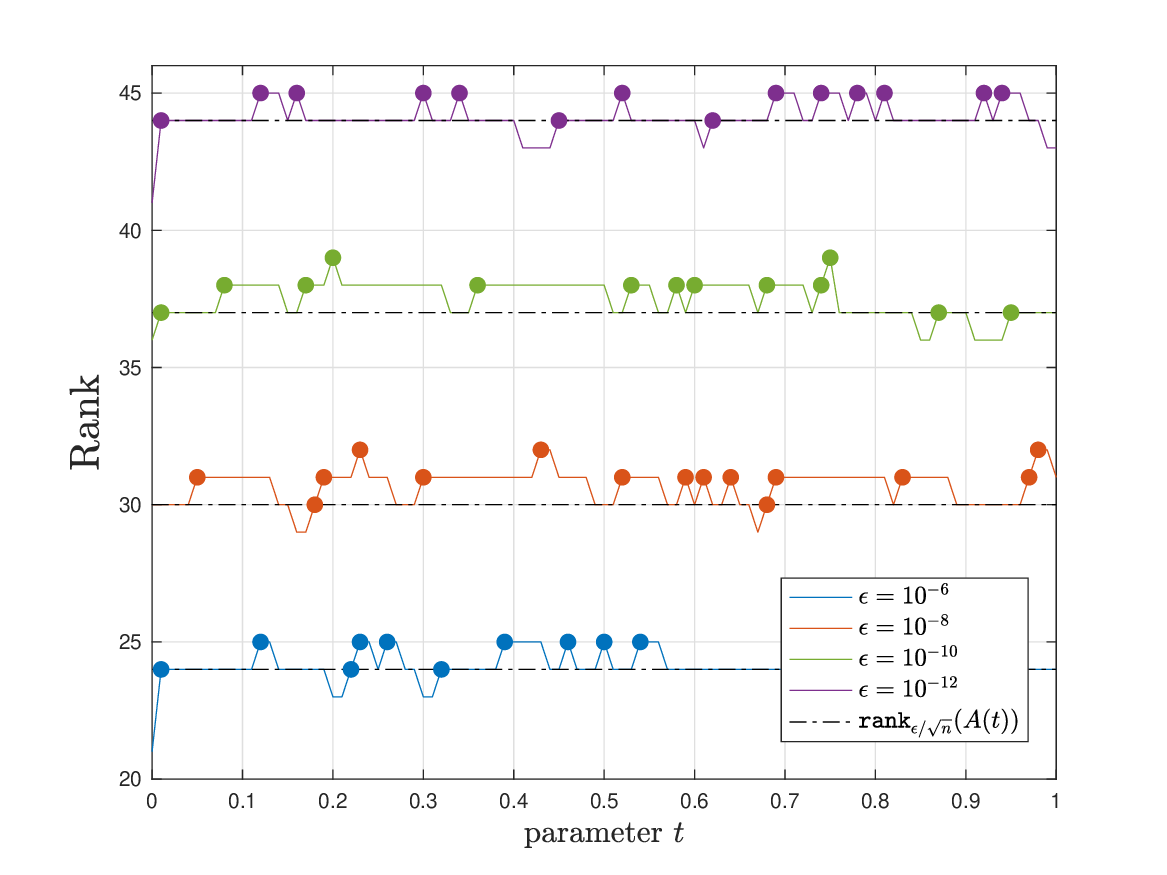}} 
\centering 
\caption{Testing tolerance parameter $\epsilon$ for \acur{} and \qcur{} using the synthetic problem \eqref{eq:synprob}.}
\label{fig:TolPlot}
\end{figure}

The results are depicted in Figure \ref{fig:TolPlot} and Table \ref{table:tol}. In Figures \ref{subfig:TolM1} and \ref{subfig:TolM2}, we run \acur{} with varying tolerance $\epsilon$ on the synthetic problem \eqref{eq:synprob}. We set the oversampling parameter $p$ to $5$ and the error sample size $s$ to $5$ in this experiment. We first observe in Figure \ref{subfig:TolM1} that the relative error for \acur{} meets the tolerance $\epsilon$ within a small constant factor. This is expected, given that the randomized methods in the algorithm typically yield errors within a small constant factor of the optimal solution; therefore one can simply take a smaller $\epsilon$ to ensure the actual error is bounded by the desired accuracy. Furthermore, the performance of \acur{} is reflected in Figure \ref{subfig:TolM2} where the algorithm attempts to efficiently find the low-rank approximation of the parameter-dependent matrix by matching the $(\epsilon/\sqrt{n})$-rank at each parameter value. Regarding the algorithm's complexity, as shown in Table \ref{table:tol}, we observe that $h_1$ and $h_2$ increases as the tolerance $\epsilon$ decreases. This increase is anticipated since a higher accuracy requirement generally equates to heavier computations.

In Figures \ref{subfig:TolA1} and \ref{subfig:TolA2}, we run \qcur{} with varying tolerance $\epsilon$ on the synthetic problem \eqref{eq:synprob}. We set the oversampling parameter $p$ to $5$ and the buffer space $b$ to $5$. We notice in Figure \ref{subfig:TolA1} that \qcur{} performs well by meeting the tolerance $\epsilon$ within a small constant factor. This performance is also reflected in Figure \ref{subfig:TolA2}, where the algorithm tries to match the $(\epsilon/\sqrt{n})$-rank at each parameter value. However, as shown in the numerical examples that follow, \qcur{} must be used with caution because as the problem becomes more difficult, the error may continue to grow as we iterate through the parameter values.

\begin{table}[!ht]
\footnotesize
\centering
\caption{The number of instances of heavier computations performed out of $100$ parameter values $t$ by \acur{} in the experiment corresponding to Figure \ref{fig:TolPlot}.}
\label{table:tol}
\begin{tabular}{|p{7cm}|p{3cm}|p{3cm}|p{3cm}|p{3cm}|}
 \hline
\multicolumn{1}{|c||}{Tolerance $\epsilon$} & \multicolumn{1}{|c|}{$10^{-6}$} & \multicolumn{1}{|c|}{$10^{-8}$} & \multicolumn{1}{|c|}{$10^{-10}$} & \multicolumn{1}{|c|}{$10^{-12}$}\\
 \hline
\multicolumn{1}{|c||}{Minor mod. only $h_1$} & \multicolumn{1}{|r|}{$6$}  & \multicolumn{1}{|r|}{$8$} & \multicolumn{1}{|r|}{$10$} & \multicolumn{1}{|r|}{$11$}   \\
\multicolumn{1}{|c||}{Recomput. Indices $h_2$} & \multicolumn{1}{|r|}{$0$}  & \multicolumn{1}{|r|}{$0$} & \multicolumn{1}{|r|}{$0$} & \multicolumn{1}{|r|}{$1$}  \\
 \hline
\end{tabular}
\end{table}

\subsection{Varying the oversampling parameter $p$} \label{subsec:RoleP}
In this part of the section, we test the role that the oversampling parameter $p$ plays in our algorithm. As shown in \cite{pn24cur}, oversampling can help with the accuracy and the numerical stability of the CUR decomposition. We use a problem based on a PDE \cite{CarrelGanderVandereycken2023,KressnerTobler2011} whose aim is to solve a heat equation simultaneously with several heat coefficients. This problem is also known as the parametric cookie problem. In matrix form, the problem is given by
\begin{equation} \label{eq:PCprob}
    \dot{A}(t) = -B_0A(t)-B_1A(t)C +b\mathbf{1}^T,\,\, A(t_0) = A_0 \in \R^{1580\times 101}
\end{equation} where $B_0,B_1\in \R^{1580\times 1580}$ are the mass and stiffness matrices, respectively, $b\in \R^{1580}$ is the discretized inhomogeneity and $C = \mathrm{diag}(0,1,...,100)$ is a diagonal matrix containing the parameter samples on its diagonal. We follow \cite{KressnerLam2023} by setting $t_0 = -0.01$ and $A_0$ as the zero matrix and compute the solution at $t=0$ exactly, after which we discretize the interval $[0,0.9]$ with $101$ uniform points and obtain $A(t)$ using the $\mathtt{ode45}$ command in MATLAB with tolerance parameters \{`RelTol',$1e-12$,`AbsTol',$1e-12$\} on each subinterval.

\begin{figure}[!ht]
\subfloat[\acur{} with $\epsilon = 10^{-12}$, $s = 1$]{\label{subfig:OSM1}\includegraphics[scale = 0.33]{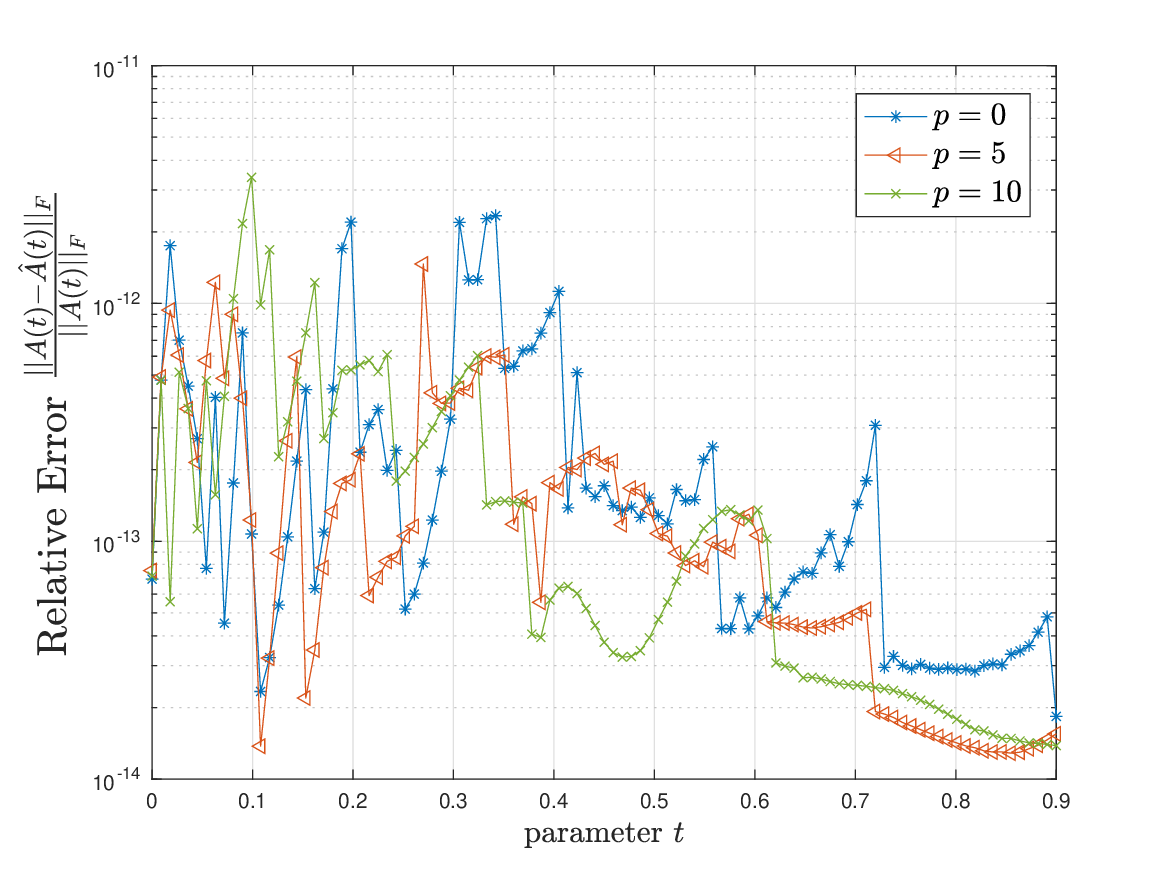}} \hspace{-0.5cm}
\subfloat[Rank changes for \acur{}]{\label{subfig:OSM2}\includegraphics[scale = 0.33]{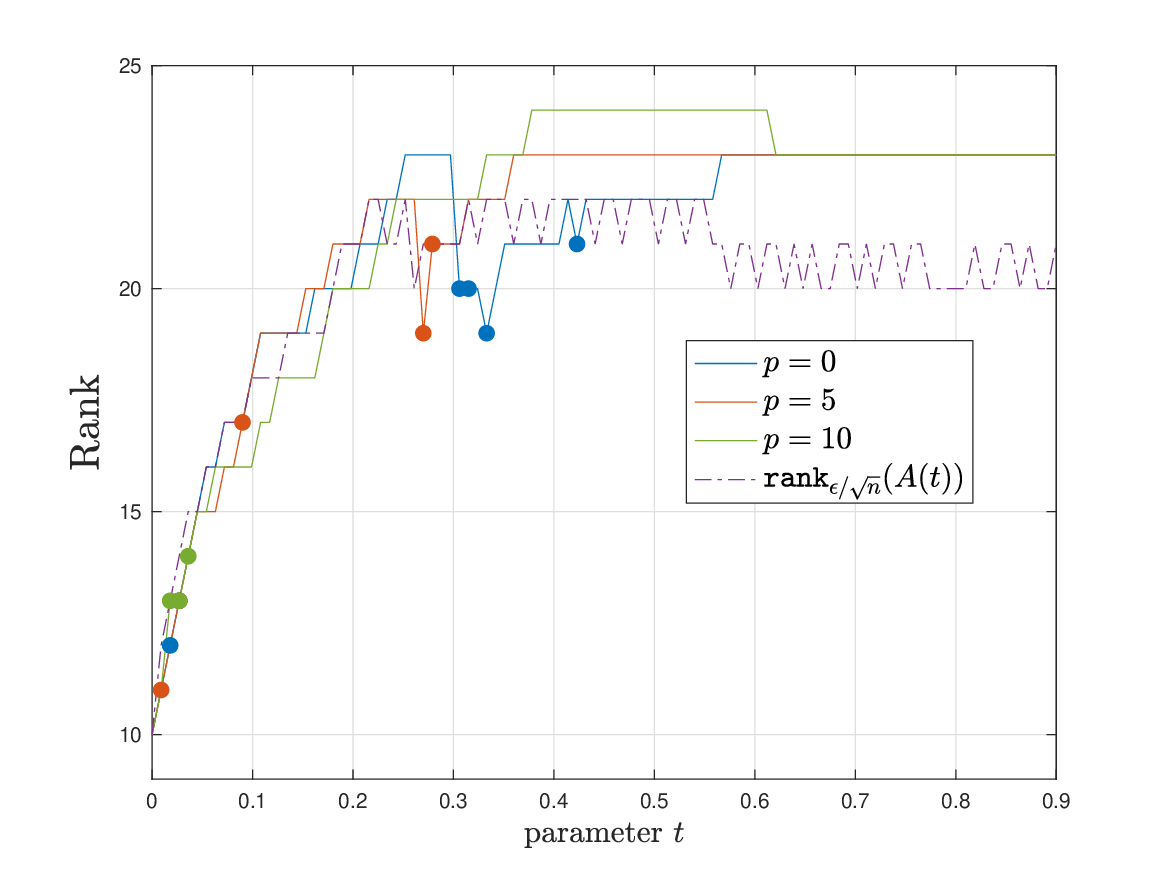}} \vspace{-0.4cm} \\
\subfloat[\qcur{} with $\epsilon = 10^{-12}$, $b = 5$]{\label{subfig:OSA1}\includegraphics[scale = 0.33]{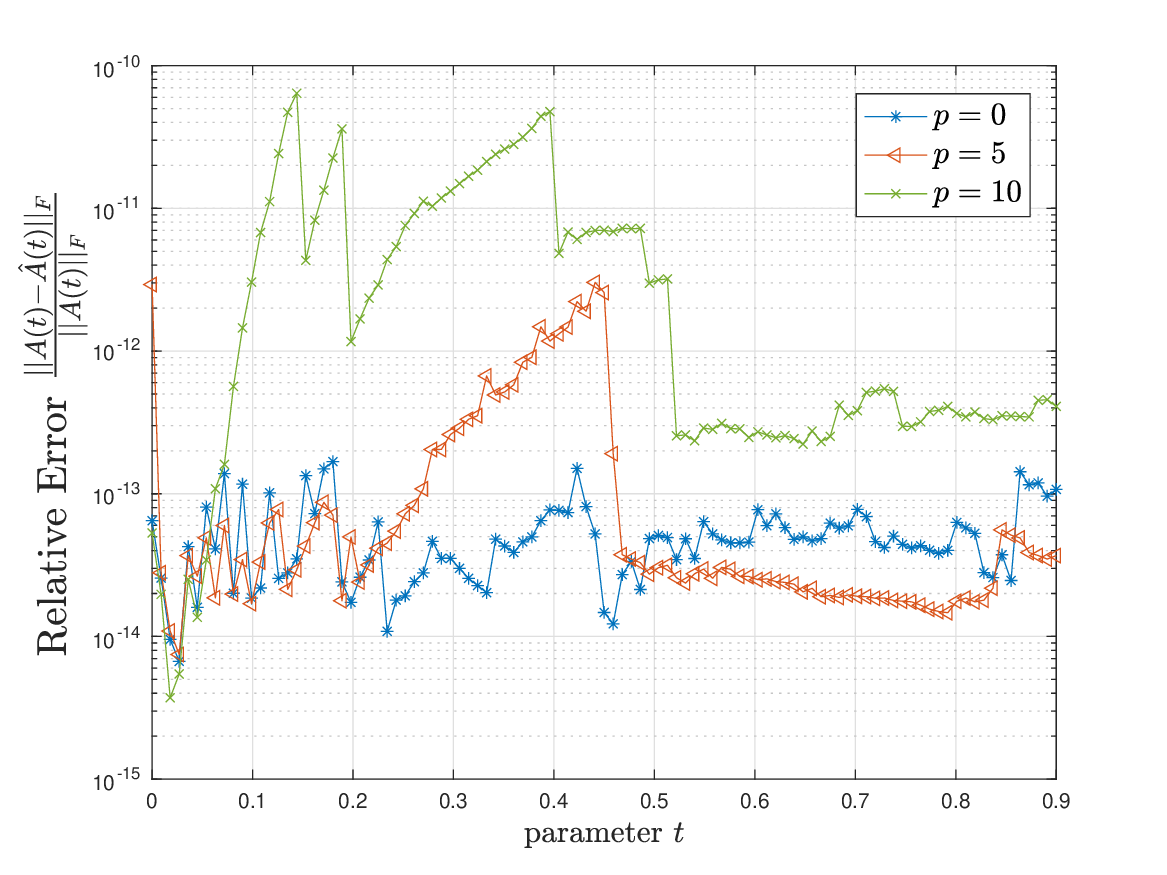}} \hspace{-0.5cm}
\subfloat[Rank changes for \qcur{}]{\label{subfig:OSA2}\includegraphics[scale = 0.33]{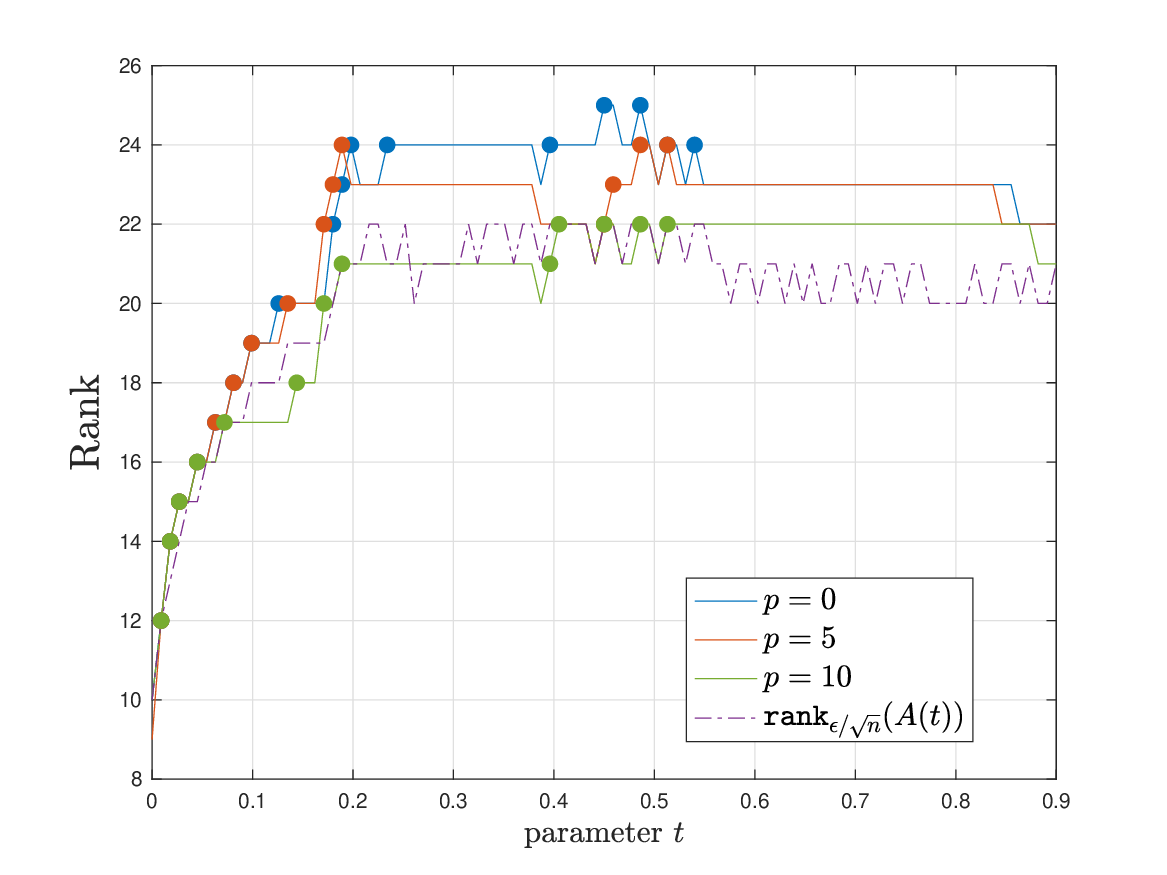}} 
\centering 
\caption{Testing the oversampling parameter $p$ for \acur{} and \qcur{} using the parametric cookie problem \eqref{eq:PCprob}.}
\label{fig:OSPlot}
\end{figure}

The results are presented in Figure \ref{fig:OSPlot} and Table \ref{table:oversample}. In Figures \ref{subfig:OSM1} and \ref{subfig:OSM2}, \acur{} was used with varying oversampling parameter $p$ on the parametric cookie problem \eqref{eq:PCprob}.
We set the tolerance $\epsilon$ to $10^{-12}$ and the error sample size $s$ to $1$ in this experiment. In Figure \ref{subfig:OSM1}, the algorithm meets the tolerance level up to a modest constant factor while approximately matching the $(\epsilon/\sqrt{n})$-rank in Figure \ref{subfig:OSM2}. The benefits of oversampling is demonstrated in Table \ref{table:oversample} where the number of instances the algorithm recomputes the indices (which forms the dominant cost), $h_2$, decreases as the oversampling parameter $p$ increases. This demonstrates that oversampling assists \acur{} in improving its complexity by reducing the frequency of heavier computations.

In Figures \ref{subfig:OSA1} and \ref{subfig:OSA2}, we use \qcur{}. We vary the oversampling parameter $p$ on the parametric cookie problem \eqref{eq:PCprob} with the tolerance $\epsilon$ equal to $10^{-12}$ and the buffer size $b$ set to $5$. In Figure \ref{subfig:OSA1}, the algorithm fails to meet the tolerance $\epsilon = 10^{-12}$ for certain parameter values, but still meets the tolerance within 2 orders of magnitude. By oversampling we are increasing the singular values of the core matrix $U_i = A(t_i)(I,J)\in \R^{(r_i+b+p)\times (r_i+b)}$; see Section \ref{subsec:cheapalg}. This, in turn, should advocate a rank increase in \qcur{}, which the algorithm uses to reduce the error. However, as the algorithm is based on heuristics to prioritize speed, the impact that the parameters have on the algorithm is rather complicated and inconclusive.\footnote{Running \qcur{} on the parametric cookie problem with the same parameter values as in Figure \ref{subfig:OSA1} several times yield different results with varying performance between different oversampling values, making the impact that oversampling has for the accuracy inconclusive.} Nevertheless, for the accuracy and stability of the CUR decomposition, we recommend at least a little bit of oversampling; see \cite{pn24cur}.

\begin{table}[!ht]
\footnotesize
\centering
\caption{The number of instances of heavier computations performed out of $100$ parameter values $t$ by \acur{} in the experiment corresponding to Figure \ref{fig:OSPlot}.}
\label{table:oversample}
\begin{tabular}{|p{7cm}|p{3cm}|p{3cm}|p{3cm}|}
 \hline
\multicolumn{1}{|c||}{Oversampling $p$} & \multicolumn{1}{|c|}{$0$} & \multicolumn{1}{|c|}{$5$} & \multicolumn{1}{|c|}{$10$}\\
 \hline
\multicolumn{1}{|c||}{Minor mod. only $h_1$} & \multicolumn{1}{|r|}{$24$}  & \multicolumn{1}{|r|}{$18$} & \multicolumn{1}{|r|}{$5$}  \\
\multicolumn{1}{|c||}{Recomput. Indices $h_2$} & \multicolumn{1}{|r|}{$6$}  & \multicolumn{1}{|r|}{$5$} & \multicolumn{1}{|r|}{$3$}  \\
 \hline
\end{tabular}
\end{table}
\subsection{Varying the error sample size $s$ and the buffer size $b$} \label{subsec:RoleS}
Here we test the role that the error sample size $s$ plays in \acur{} and the buffer size $b$ in \qcur{}. We use the discrete Schr\"odinger equation in imaginary time for this experiment. The example is taken from \cite{CerutiLubich2022} and is given by
\begin{equation} \label{eq:SchroEq}
    \dot{A}(t) = \frac{1}{2}(D\, A(t)+A(t)\, D)-V_{\cos}\, A(t)\, V_{\cos}, \,\, A(0) = A_0,
\end{equation} where $D = \mathrm{tridiag}(-1,2,-1)$ is the discrete 1D Laplacian and $V_{\cos}$ is the diagonal matrix with entries $1-\cos(2j\pi/n)$ for $j = -n/2,...,n/2-1$. We take $n = 512$ and the initial condition $A_0$ is a randomly generated matrix with singular values $10^{-i}$ for $i = 1,...,512$. We obtain $A(t)$ by discretizing the interval $[0,0.1]$ with $101$ uniform points and using the $\mathtt{ode45}$ command in MATLAB with tolerance parameters \{`RelTol',$1e-12$,`AbsTol',$1e-12$\} on each subinterval.

\begin{figure}[!ht]
\subfloat[\acur{} with $\epsilon = 10^{-12}$, $p = 10$]{\label{subfig:SSM1}\includegraphics[scale = 0.33]{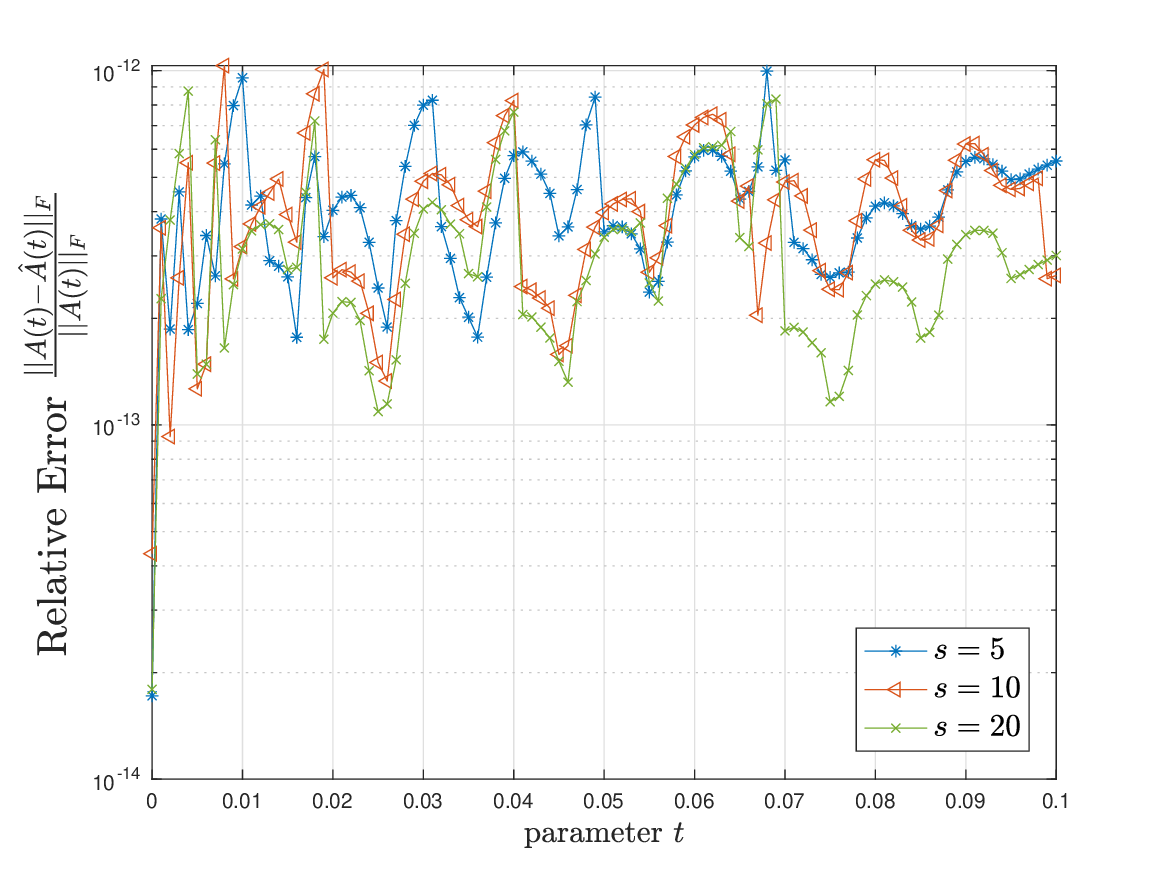}} \hspace{-0.5cm}
\subfloat[Rank changes for \acur{}]{\label{subfig:SSM2}\includegraphics[scale = 0.33]{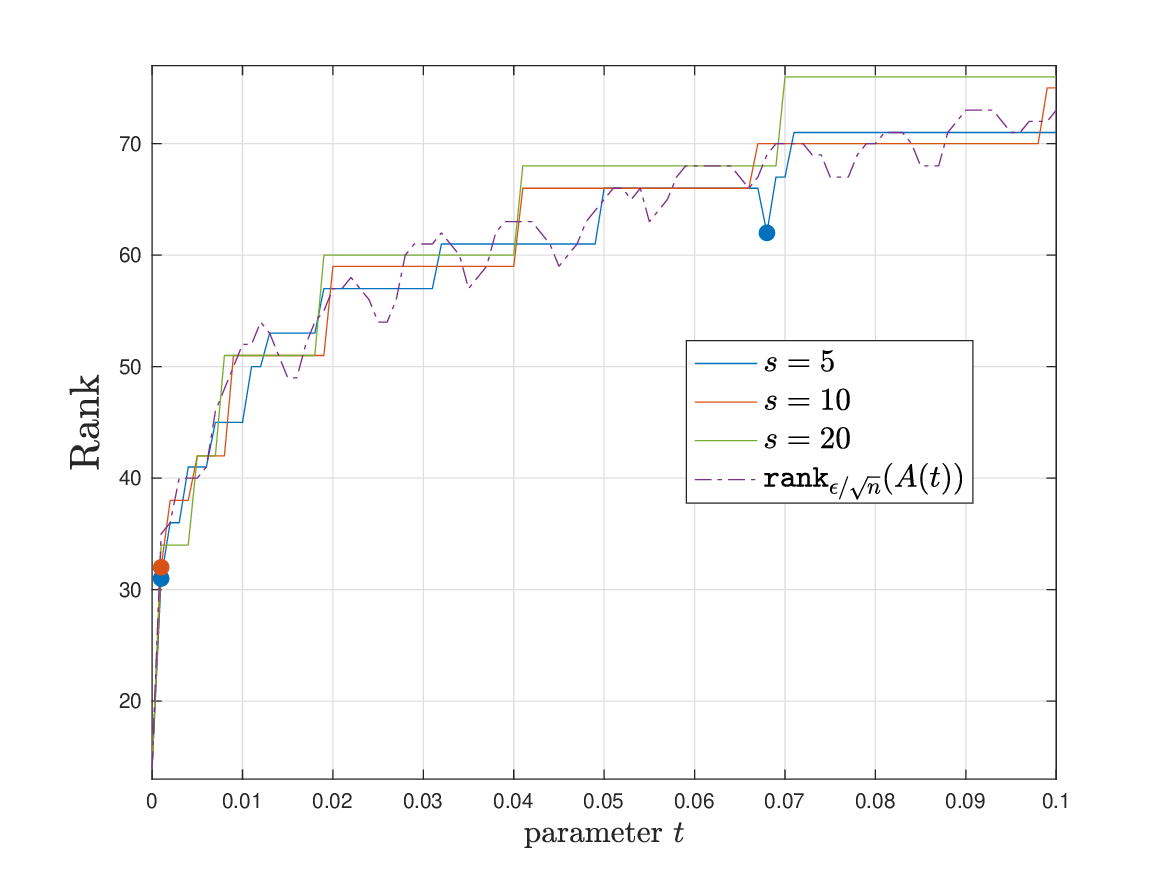}} \vspace{-0.4cm} \\
\subfloat[\qcur{} with $\epsilon = 10^{-12}$, $p = 10$]{\label{subfig:BSA1}\includegraphics[scale = 0.33]{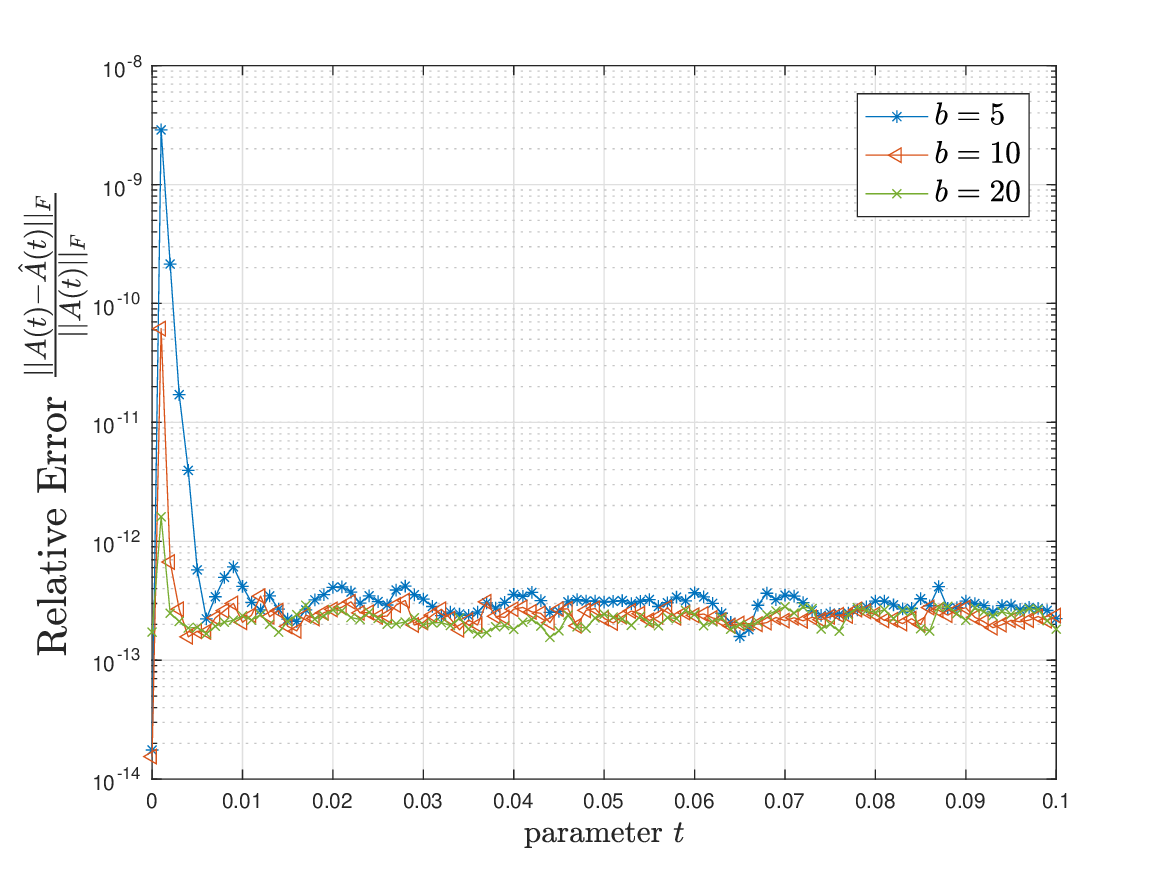}} \hspace{-0.5cm}
\subfloat[Rank changes for \qcur{}]{\label{subfig:BSA2}\includegraphics[scale = 0.33]{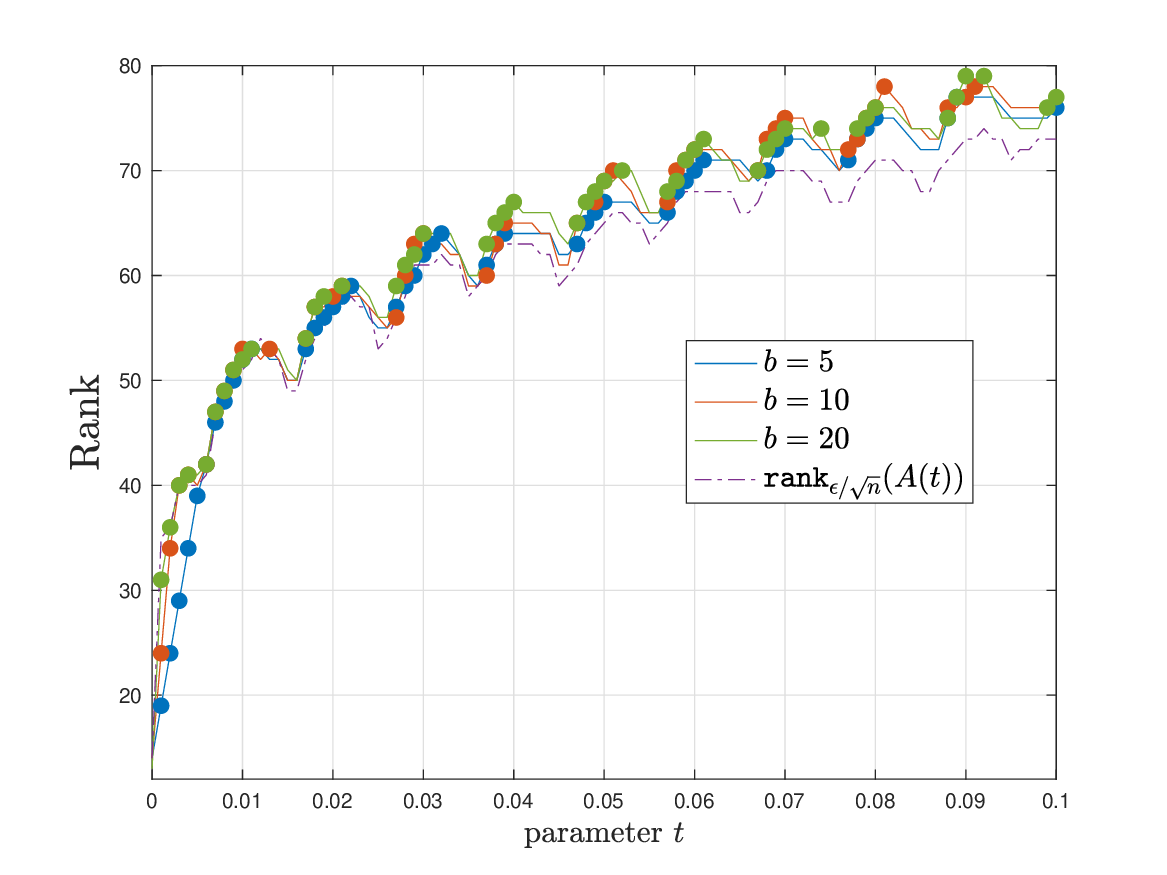}} 
\centering 
\caption{Testing the error sample size $s$ for \acur{} and the buffer size $b$ for \qcur{} using the Schr\"odinger equation \eqref{eq:SchroEq}.}
\label{fig:SPlot}
\end{figure}

The results are shown in Figure \ref{fig:SPlot} and Table \ref{table:samplesize}. In Figures \ref{subfig:SSM1} and \ref{subfig:SSM2}, we execute \acur{} with varying error sample size $s$ on the Schr\"odinger equation \eqref{eq:SchroEq}. The oversampling parameter was set to $p = 10$ and the tolerance $\epsilon = 10^{-12}$. The algorithm satisfies the tolerance $\epsilon$ up to a modest constant factor as demonstrated in Figure \ref{subfig:SSM1}, while approximately aligning with the $(\epsilon/\sqrt{n})$-rank, as depicted in Figure \ref{subfig:SSM2}. As described in the role of the error sample size $s$ in Section \ref{subsec:mainalg}, it makes the minor modifications in \acur{} more effective. This is demonstrated in Table \ref{table:samplesize} where the number of times the algorithm recomputed the indices, $h_2$, decreases as the error sample size $s$ increases. The experiment demonstrates that with an increase in $s$, the algorithm becomes better at lowering the relative error using the low-cost minor modifications only. Therefore, the error sample size $s$ should be set sensibly to allow the algorithm to resolve the inaccuracy using minor modifications only so that the algorithm does not recompute indices until it becomes necessary. Note that higher $s$ increases the complexity of \acur{}, so $s$ should not be too large.

In Figures \ref{subfig:BSA1} and \ref{subfig:BSA2}, we execute \qcur{}. We vary the buffer size $b$ on the Schr\"odinger equation \ref{eq:SchroEq} with the tolerance set to $\epsilon = 10^{-12}$ and the oversampling parameter set to $p = 10$. In Figure \ref{subfig:BSA1}, the algorithm performs well by meeting the tolerance except for the large spike happening in the first few parameter values. The cause of this spike is explained in Figure \ref{subfig:BSA2} where the problem experiences a substantial rank increase in the initial parameter values. Since the algorithm is only able to make a maximum rank increase of $b$ at each iteration, the algorithm struggles to keep up with the large rank increase in the initial parameter values. This observation is evident in Figure \ref{subfig:BSA1}, where the initial spike in the graph diminishes as $b$ increases. Therefore, if we expect a large rank increase in the parameter-dependent matrix, the buffer size $b$ should be set to a higher value for those parameter values.

\begin{table}[!ht]
\footnotesize
\centering
\caption{The number of instances of heavier computations performed out of $100$ parameter values $t$ by \acur{} in the experiment corresponding to Figure \ref{fig:SPlot}.}
\label{table:samplesize}
\begin{tabular}{|p{7cm}|p{3cm}|p{3cm}|p{3cm}|}
 \hline
\multicolumn{1}{|c||}{Error sample size $s$} & \multicolumn{1}{|c|}{$5$} & \multicolumn{1}{|c|}{$10$} & \multicolumn{1}{|c|}{$20$}\\
 \hline
\multicolumn{1}{|c||}{Minor mod. only $h_1$} & \multicolumn{1}{|r|}{$12$}  & \multicolumn{1}{|r|}{$8$} & \multicolumn{1}{|r|}{$6$}  \\
\multicolumn{1}{|c||}{Recomput. Indices $h_2$} & \multicolumn{1}{|r|}{$2$}  & \multicolumn{1}{|r|}{$1$} & \multicolumn{1}{|r|}{$0$}  \\
 \hline
\end{tabular}
\end{table}

\subsection{Adversarial example for \qcur{}} \label{subsec:adversarial}
In this section, we create an adversarial example for \qcur{} to show that the algorithm can fail. We propose the following block diagonal matrix
    \begin{equation}
        A(t) = \begin{bmatrix}
            A_1 & 0 & 0\\ 0 & 0 & c(t) t A_2 
        \end{bmatrix} \in \R^{300\times 100}, t\in [0,1]
    \end{equation} where $A_1 = \mathtt{randn}(100,20)$, $A_2 = \mathtt{randn}(200,10)$ and $c(t) = 10^{-5+10t}$ using the $\mathtt{randn}$ command in MATLAB. The parameter-dependent matrix $A(t)$ starts with a large component in the $A_1$ block, since $c(t)tA_2$ is the zero matrix at $t=0$. As $t$ increases, the dominant block becomes $c(t)t A_2$. 
    
\begin{figure}[!ht]
\subfloat[\acur{} with $\epsilon = 10^{-4}$]{\label{subfig:AM1}\includegraphics[scale = 0.33]{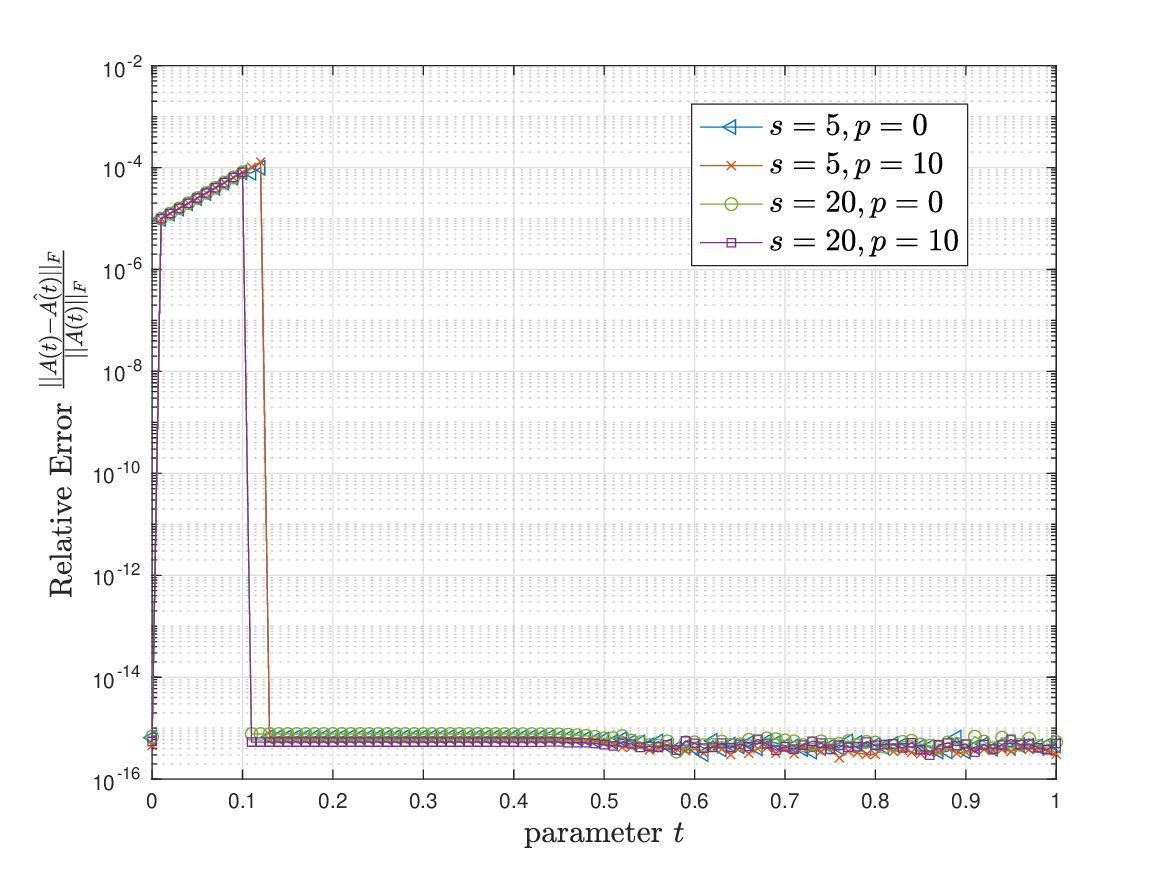}} \hspace{-0.5cm}
\subfloat[\qcur{} with $\epsilon = 10^{-4}$]{\label{subfig:AA1}\includegraphics[scale = 0.33]{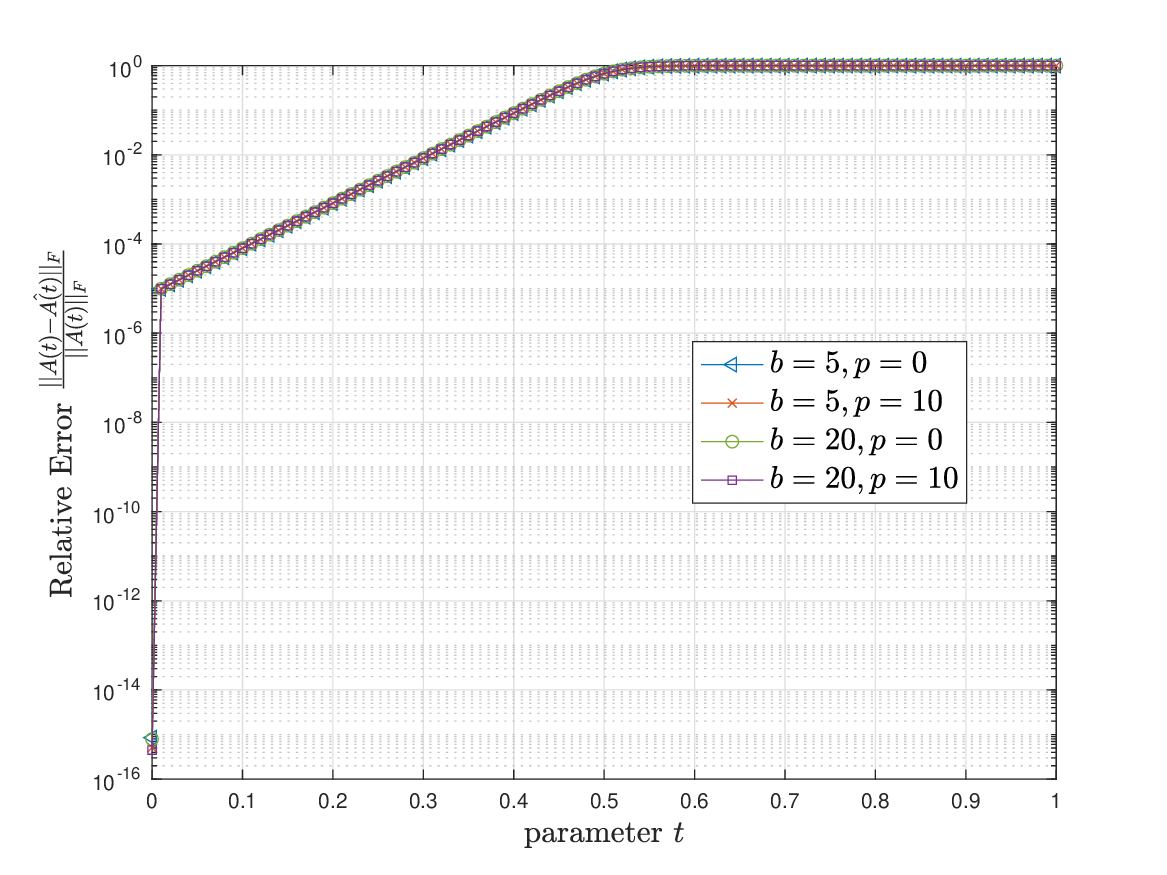}} 
\centering 
\caption{Adversarial example for \qcur{} with tolerance $\epsilon = 10^{-4}$. The adversarial example causes \qcur{} to fail, whereas \acur{} remains successful.}
\label{fig:AdPlot}
\end{figure}

    The results are presented in Figure \ref{fig:AdPlot}. We execute \qcur{} in Figure \ref{subfig:AA1} with tolerance $\epsilon = 10^{-4}$ and various values of buffer size $b$ and oversampling parameter $p$. \qcur{} fails as the error continues to grow as we iterate through the parameter values. This can be explained by the nature of the algorithm. Since it only considers a submatrix of the original matrix for each parameter value, the $c(t)tA_2$ block remains hidden as $t$ increases. Consequently, the algorithm is never able to capture the $c(t)tA_2$ block, making the approximation poor.\footnote{This type of counterexample is present in essentially all algorithms that do not view the entire matrix. For example, the algorithm in \cite{DonelloPalkarNaderiDelReyFernandezBabaee2023} may suffer from the same adversarial problem as it does not view the entire matrix at each parameter value. See \cite[\S~17.5]{martinsson2019fast} for a related discussion.} On the other hand, as illustrated in Figure \ref{subfig:AM1}, \acur{} remains successful for the adversarial problem. \acur{} is able to satisfy the tolerance of $10^{-4}$, and when the $c(t)tA_2$ block starts to become dominant, minor modifications or recomputation of indices is employed to capture the $c(t)tA_2$ block, making the algorithm successful.
    
\subsection{Comparison against other methods}
In this section, we test the data-driven aspect of \acur{} and \qcur{} and compare their speed against other methods. Specifically, we benchmark \acur{} and \qcur{} against the randomized SVD and the generalized Nystr\"om method from \cite{KressnerLam2023}, as well as an algorithm in the dynamical low-rank approximation literature proposed by Lubich and Oseledets \cite{LubichOseledets2014}.

\subsubsection{Rank-adaptivity test}
To test the rank-adaptive nature of \acur{} and \qcur{}, we use the following synthetic problem. Let $A\in \R^{1000\times 300}$ be a low-rank matrix with $\rk(A) = 50$, Haar-distributed left and right singular vectors, and singular values that decay geometrically from $1$ to $10^{-12}$. We define the parameter-dependent matrix $A(t)$ as
\begin{equation}
    A(i) = A + 10^{-4}\sum_{j = 2}^{i} G_1^{(j)} G_2^{(j)}
\end{equation} where $i \in \{1,2,...,10\}$, and $G_1^{(j)}\in \R^{1000\times 4}$, $G_2^{(j)}\in \R^{4\times 300}$ are Gaussian matrices with i.i.d. entries $\mathcal{N}(0,1/1000)$ and $\mathcal{N}(0,1/300)$ respectively. This experiment starts with a rank-$50$ matrix and incrementally adds rank-$4$ perturbations of size $10^{-4}$. 

For \acur{} and \qcur{}, the input parameters were set as follows: $\epsilon = 10^{-9}$ for tolerance, $s,b = 5$ for the error sample size and the buffer size, and $p=10$ for the oversampling parameter. For the randomized SVD, generalized Nystr\"om method and the Lubich-Oseledets algorithm, the target rank was set to $50$, corresponding to the rank of $A$. The second sketch size of the generalized Nystr\"om was set to $75$. See their respective papers \cite{KressnerLam2023,LubichOseledets2014} for further details.

\begin{figure}[!ht]
\hspace{-0.4cm}
\subfloat{\includegraphics[scale = 0.35]{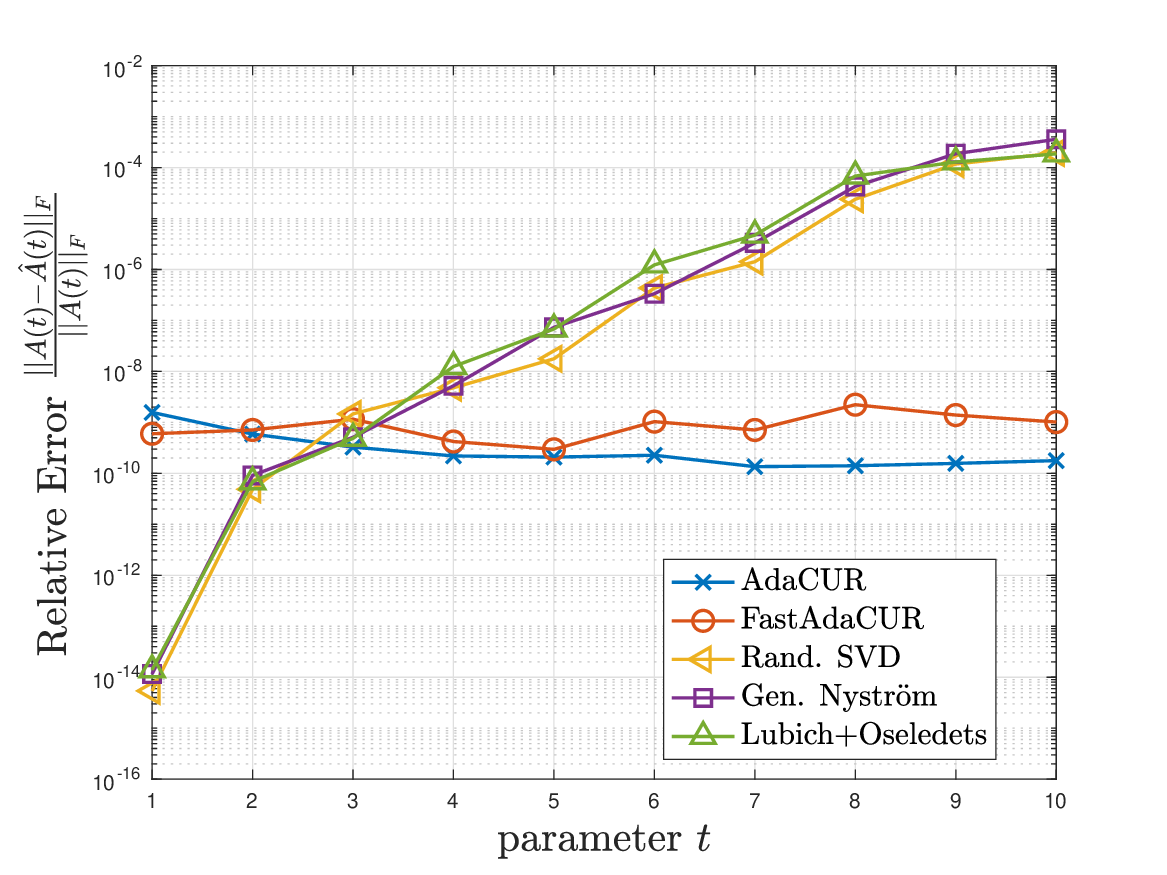}} \hspace{-0.7cm}
\subfloat{\includegraphics[scale = 0.35]{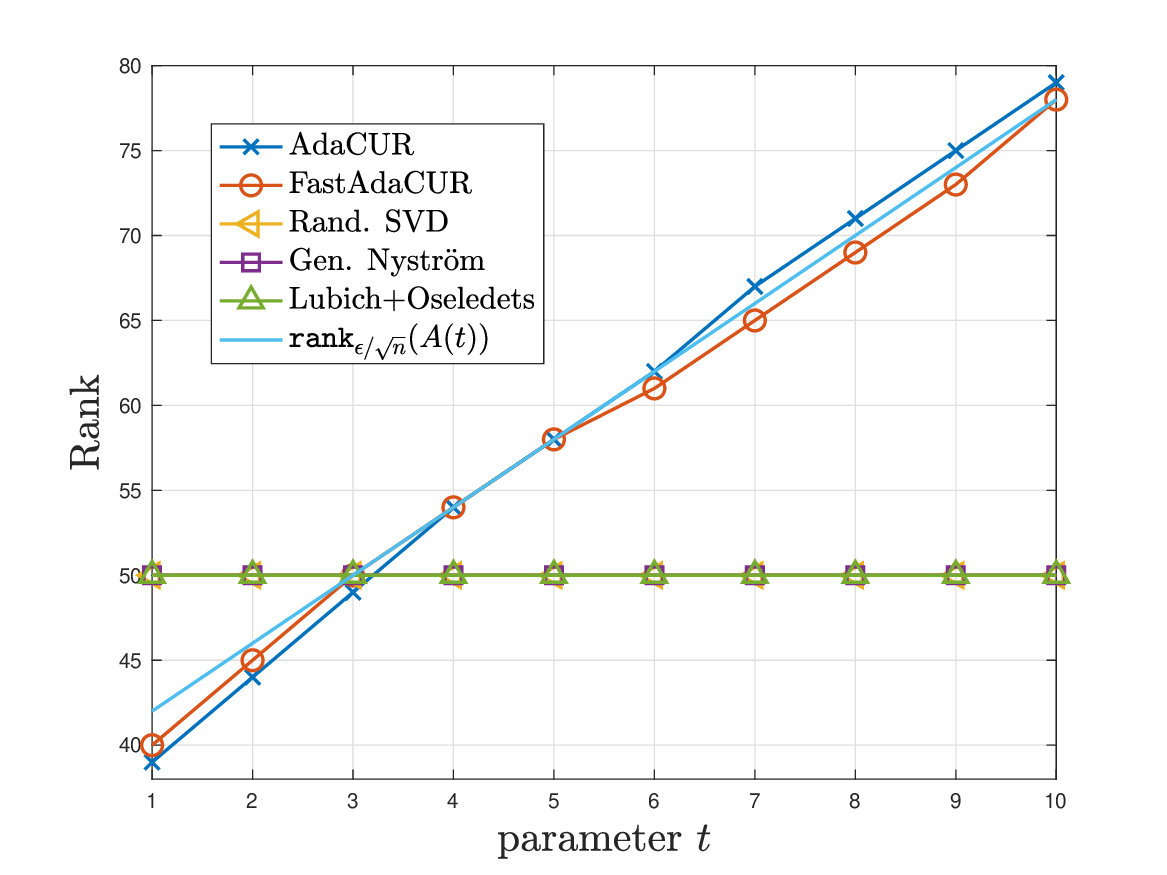}} 
\centering 
\caption{Rank-adaptivity test of \acur{} and \qcur{} against the three existing algorithms: randomized SVD, generalized Nystr\"om method \cite{KressnerLam2023} and the algorithm by Lubich and Osedelets \cite{LubichOseledets2014}.}
\label{fig:dataPlot}
\end{figure}

The results are shown in Figure~\ref{fig:dataPlot}. We observe that fixed-rank methods--randomized SVD, generalized Nystr\"om, and the Lubich-Oseledets algorithm--perform poorly as parameter value changes, failing to adapt to evolving data. In contrast, \acur{} and \qcur{} effectively adjust the rank as needed to try maintain the error below the specified tolerance. This is guaranteed for \acur{} with high probability due to its error control mechanism. This experiment highlights the data-driven nature of \acur{} and \qcur{}, which dynamically adapt to the data, unlike the fixed-rank approaches. While \acur{} and \qcur{} are effective rank-adaptive algorithms, they are not the only methods for addressing rank-adaptivity. Notably, rank-adaptivity has been explored in the dynamical low-rank approximation literature \cite{CerutiKuschLubich2022,HauckSchnake2023,HesthavenPagliantiniRipamonti2022,HochbruckNeherSchrammer2023}. This experiment highlights the simplicity and effectiveness of \acur{} and \qcur{} in adapting to rank changes.

\subsubsection{Speed test}\label{subsec:Speed}
When the target rank $r$ becomes sufficiently large, we anticipate that our algorithms, which run with $\bigO(T_{A(t)})$ and $\bigO((m+n)r^2)$ complexity, will outperform many existing methods, which run with $\bigO(r T_{A(t)})$ complexity. We demonstrate this through experiments using the following artificial problem. Let $A$ be a low-rank matrix with $\rk(A) = 500$ given by $A = U\Sigma V^T \in \R^{50000\times 5000}$ where $U\in \R^{50000\times 500}$ and $V\in \R^{5000 \times 500}$ are Haar-distributed orthonormal matrices and $\Sigma \in \R^{500\times 500}$ is a diagonal matrix with entries that decay geometrically from $1$ to $10^{-8}$. The parameter-dependent matrix $A(t)$ is then defined by
\begin{equation}
    A(i) = A + \delta \sum_{j=1}^{i-1} X_j
\end{equation} where $i\in \{1,2,...,101\}$, $\delta = 10^{-12}$ and $X_j$'s are sparse random matrices generated using $\mathtt{sprandn(50000,5000,10^{-5})}$ command in MATLAB. The parameter-dependent matrix $A(t)$ starts as a low-rank matrix at $i=1$, with sparse noise introduced at discrete time intervals. \acur{} and \qcur's input parameters were $\epsilon = 10^{-6}$ for tolerance, $s,b = 10$ for the error sample size and the buffer size, and the oversampling parameter was set to $p = 10$. The target rank for the randomized SVD, generalized Nystr\"om method and the algorithm by Lubich and Osedelets was set to the rank of $A$, which is $500$. The second sketch size of the generalized Nystr\"om was set to $750$. See their respective papers for details \cite{KressnerLam2023,LubichOseledets2014}.

\begin{figure}[!ht]
\subfloat{\includegraphics[scale = 0.33]{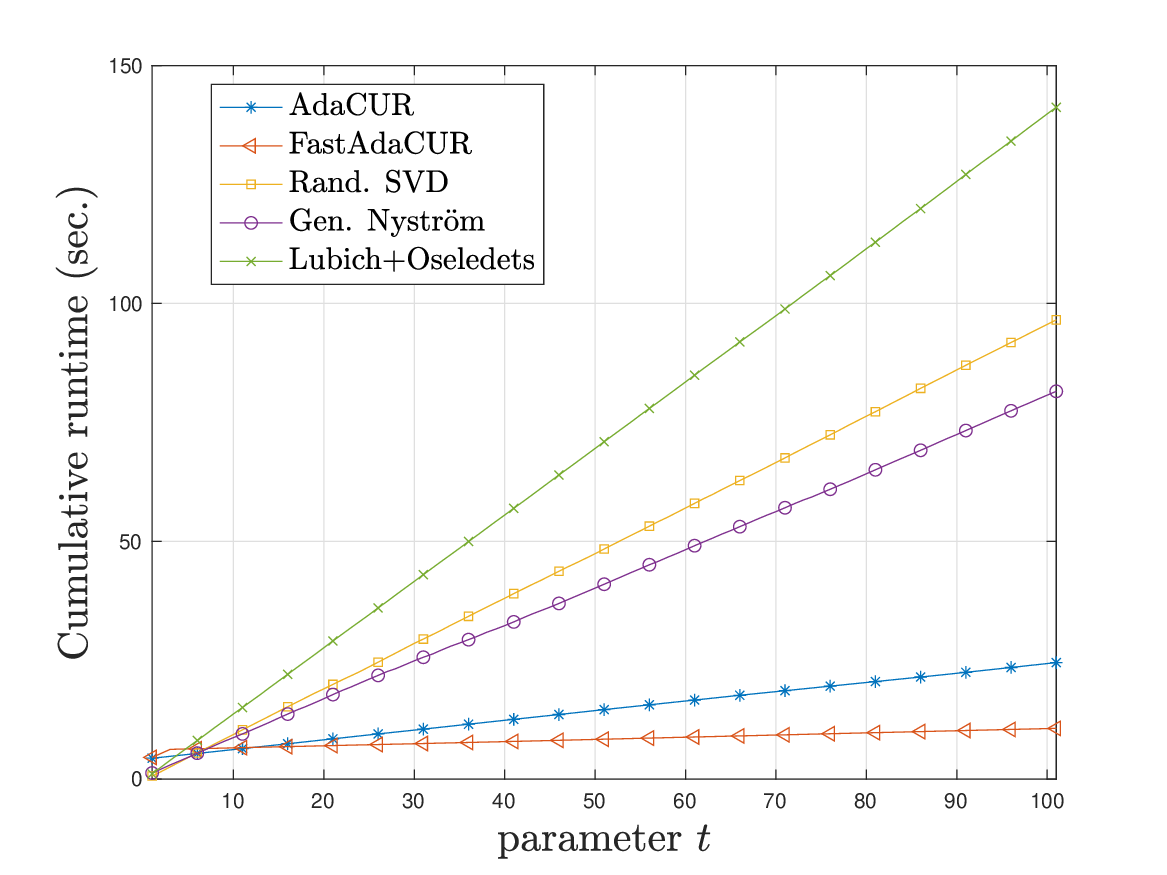}} \hspace{-0.5cm}
\subfloat{\includegraphics[scale = 0.33]{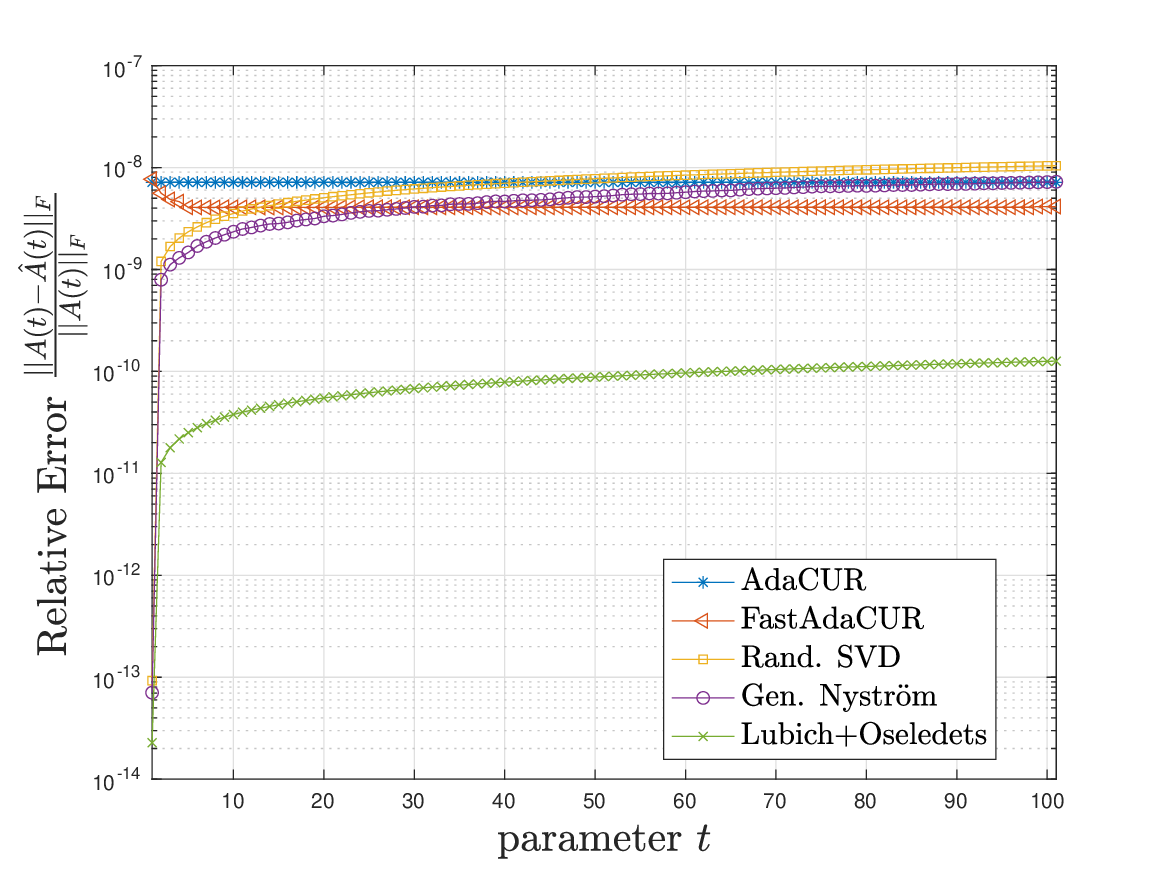}} 
\centering 
\caption{Runtime test for \acur{} and \qcur{} against the three existing algorithms: randomized SVD, generalized Nystr\"om method \cite{KressnerLam2023} and the algorithm by Lubich and Osedelets \cite{LubichOseledets2014}.}
\label{fig:SpeedPlot}
\end{figure}

The results are illustrated in Figure~\ref{fig:SpeedPlot}. For \acur{} and \qcur{}, we notice a slight rise at the beginning, stemming from the initial computation of indices. However, the low-rank approximation of the subsequent parameter values is computed very quickly, making the slope of the cumulative runtime graph relatively flat as we iterate through the parameter values. The graph is notably flat for \qcur{}, as it runs at most linearly in $m$ and $n$ at each iteration. On the other hand, the three existing methods display a steeper slope due to the higher computational cost associated with each parameter value. This higher cost stems from the increased number of matrix-vector products with $A(t)$ or its transpose. \acur{} is approximately $3$--$5$ times faster than the three existing algorithms, while \qcur{} is approximately $7$--$13$ times faster than the three existing algorithms.

\section{Conclusion}
In this work, we devised two efficient rank-adaptive algorithms for computing the low-rank approximation of parameter-dependent matrices: \acur{} and \qcur{}. The key idea behind these algorithms is to try to reuse the row and the column indices from the previous iterate as much as possible. \acur{} comes with many favourable properties such as rank-adaptivity, error-control and a typical complexity of $\bigO(T_{A(t)})$, while \qcur{}, which is also rank-adaptive, is faster with a complexity that is at most linear in $m$ and $n$, but lacks error control, making it susceptible to adversarial problems. Nonetheless, \qcur{} should work on easier problems such as parameter-dependent matrices that are coherent at each parameter value.

The challenge in \qcur{} is that we require an efficient way to detect large errors without viewing the entire matrix. Adversarial examples such as the one in Section \ref{subsec:adversarial} always exist if we do not view the entire matrix. However, a probabilistic method that allows some sort of error control with high probability would be beneficial for \qcur{}. This is left for future work.

\bibliographystyle{siamplain}
\bibliography{references}

\end{document}